\documentclass[a4paper,twoside,10pt]{article}
\usepackage{amssymb}
\usepackage{amsmath}
\usepackage{amsfonts}
\usepackage{color}
\usepackage{mathrsfs}
\usepackage{amsmath,amsfonts,amssymb,amsthm}
\usepackage{enumitem}
\usepackage[bookmarks]{hyperref}

\usepackage[raggedright]{titlesec}
\usepackage{tocbibind}
\usepackage{geometry}
\numberwithin{equation}{section}

\newtheorem{thm}{Theorem}[section]
\newtheorem{cj}{Conjecture}[section]
\newtheorem{defi}[thm]{Definition}
\newtheorem{lem}[thm]{Lemma}
\newtheorem{cor}[thm]{Corollary}
\newtheorem{rem}[thm]{Remark}
\newtheorem{prop}[thm]{Proposition}
\newtheorem{nota}[thm]{Notation}

\def\bb{\mathbb}
\def\ca{\mathcal}

\def\e{\epsilon}
\def\cL{{\mathcal L}}
\def\E{\mathbb{E}}
\def\cW{{\mathcal W}}

\def\scr{\mathscr}
\def\bbd{\mathbf}

\def\dirac{\boldsymbol{\delta}}

\def\bthm{\begin{thm}\def\ethm{\end{thm}}}
\def\brem{\begin{rem}\def\erem{\end{rem}}}

\def\bdefi{\begin{defi}\sl{}\def\edefi{\end{defi}}}
\def\blem{\begin{lem}\def\elem{\end{lem}}}

\def\bprop{\begin{prop}\def\eprop{\end{prop}}}

\def\bitm{\begin{itemize}} \def\eitm{\end{itemize}}
\def\benu{\begin{enumerate}} \def\eenu{\end{enumerate}}

\def\bpf{\begin{proof}}\def\epf{\end{proof}}
\def\beq{\begin{equation}}\def\eeq{\end{equation}}
\def\beqs{\begin{eqnarray}}\def\eeqs{\end{eqnarray}}
\def\beqsnl{\begin{eqnarray*}}\def\eeqsnl{\end{eqnarray*}}

\usepackage{fancyhdr}
\begin{document}
\title{Uniqueness and propagation  of  chaos for the Boltzmann  equation  with  moderately
soft potentials}

\author{Liping Xu}

\maketitle
\begin{abstract}
We prove a strong/weak stability estimate for the 3D homogeneous Boltzmann equation with
moderately soft potentials ($\gamma\in(-1,0)$) using the Wasserstein distance with quadratic
cost. This in particular implies uniqueness in the class of all weak solutions,
assuming only that the initial condition has a finite entropy and a finite moment of
sufficiently high order.
We also consider the Nanbu $N$-stochastic
particle system which approximates the weak solution.
We use a probabilistic coupling method and
give, under suitable assumptions on the initial condition, a rate of
convergence of the empirical measure of the particle
system to the solution of the Boltzmann equation for this singular interaction.
\end{abstract}

\section{Introduction}
\subsection{The Boltzmann equation}
We consider a 3-dimensional spatially homogeneous  Boltzmann equation, which depicts the density $f_t(v)$ of particles in a gas, moving with velocity $v\in\bb{R}^3$ at time $t\geq 0$. The density $f_t(v)$ solves  \beq\label{Bol}
\partial_tf_t(v)=\int_{\bb{R}^3}dv_{\ast}\int_{\bb{S}^2}d\sigma B(|v-v_{\ast}|,\theta)[f_{t}(v^{\prime})f_{t}(v_{\ast}^{\prime})-f_t(v)f_t(v_{\ast})],
\eeq
where
\beq\label{Bol1}
v^{\prime}=\frac{v+v_{\ast}}{2}+\frac{|v-v_{*}|}{2}\sigma,~v_{\ast}^{\prime}
=\frac{v+v_{\ast}}{2}-\frac{|v-v_{\ast}|}{2}\sigma,
\eeq
and $\theta$ is the \emph{deviation angle} defined by $\cos\theta=\frac{v-v_{\ast}}{|v-v_{\ast}|}\cdot\sigma$. The  \emph{collision Kernel} $B(|v-v_{\ast}|,\theta)\geq 0$ depends on the type of interaction between particles. It only depends on $|v-v_{\ast}|$ and on the cosine of the deviation angle $\theta$.
Conservations of mass, momentum and kinetic energy hold for reasonable solutions and we may assume without loss of generality that $\int_{\mathbb{R}^3}f_{t}(v)dv=1$ for all $t\ge0$.

\subsection{Assumptions}
We will assume that there is a measurable function $\beta:(0,\pi]\rightarrow\bb{R}_+$ such that
\beq\label{con}
\left\{\begin{array}{l} B(|v-v_{\ast}|,\theta)\sin\theta={|v-v_{\ast}|}^{\gamma}\beta(\theta),\\
\exists~0<c_0<c_1,~\forall~\theta\in(0,\pi/2), ~c_0\theta^{-1-\nu}\leq\beta(\theta)\leq c_1\theta^{-1-\nu},\\
\forall~\theta\in[\pi/2,\pi],~\beta(\theta)=0,
\end{array}
\right.\eeq
 for some $\nu\in(0,1)$, and $\gamma\in(-1,0)$ satisfying $\gamma+\nu>0$.

\vskip1mm

 The last assumption $\beta=0$ on $[\pi/2,\pi]$  is not a restriction and can be obtained by  symmetry as noted in the introduction of \cite{MR1765272}.
This assumption corresponds to  a classical physical example, inverse power laws interactions:
when particles collide by pairs due to a repulsive force proportional to $1/r^s$ for some $s>2$, assumption \eqref{con} holds with $\gamma=(s-5)/(s-1)$ and $\nu=2/(s-1)$.
Here we will focus on the case of moderately soft potentials, i.e. $s\in (3,5)$.

\subsection{Some notations}
Let  us denote by $\ca{P}(\bb{R}^3)$  the set of probability measures on $\bb{R}^3$ and by $Lip(\bb{R}^3)$ the set of  bounded  globally  Lipschitz functions $\phi: \bb{R}^3 \mapsto \bb{R}$.  When $f\in\ca{P}(\bb{R}^3)$ has a density, we also denote this density  by $f$ .
For $q>0$, we set
\[\ca{P}_q(\bb{R}^3)=\{f\in\ca{P}(\bb{R}^3): m_q(f)<\infty\} \quad\text{with}~~m_q(f):=\int_{\bb{R}^3}|v|^q f(dv).\]
We now introduce, for $\theta\in(0,\pi/2)$ and $z\in[0,\infty)$,

\beq\label{nota}
H(\theta)=\int_\theta^{\pi/2}\beta(x)dx\quad\text{and}\quad G(z)=H^{-1}(z).
\eeq
Under \eqref{con}, it is clear that $H$ is a continuous decreasing function valued in  $[0,\infty)$, so  it has an inverse function $G:[0,\infty)\mapsto(0,\pi/2)$ defined by $G(H(\theta))=\theta$  and $H(G(z))=z$. Furthermore, it is easy to verify  that there exist some constants $0<c_2<c_3$ such that for all $z>0$,
\beq\label{con1}
c_2(1+z)^{-1/\nu}\le G(z)\le c_3(1+z)^{-1/\nu},
\eeq
and we know from \cite{MR2398952} that there exists a constant $c_4>0$ such that for all $x,y\in\bb{R}_+$,
\beq\label{ineq}
\int_0^{\infty}(G(z/x)-G(z/y))^2dz\le c_4\frac{(x-y)^2}{x+y}.
\eeq

Let us now introduce the Wasserstein distance with quadratic cost on $\ca{P}_2(\bb{R}^3)$.
For $g,\tilde{g}\in\ca{P}_2(\bb{R}^3)$, let $\ca{H}(g,\tilde{g})$ be the set of probability measures on $\bb{R}^3\times\bb{R}^3$ with first marginal $g$ and second marginal $\tilde{g}$. We then set
\[ \ca{W}_2(g,\tilde{g}) =  \inf\left\{\Big(\int_{\bb{R}^3\times\bb{R}^3}|v-\tilde{v}|^2R(dv,d\tilde{v})\Big)^{1/2},~~R\in\ca{H}(g,\tilde{g})\right\}. \]
For more details on this distance, one can see \cite[Chapter 2]{MR1964483}.

\subsection{Weak solutions}
We  now introduce a suitable spherical parameterization of \eqref{Bol1} as in \cite{MR1885616}. For each $x \in\bb{R}^{3}\setminus{\{0\}}$,
we consider a vector $I(x)\in\bb{R}^3$ such that $|I(x)|=|x|$ and $I(x)\perp x$. We also set $J(x)=\frac{x}{|x|}\wedge I(x)$, where $\wedge$ is the vector product. Then the triplet $(\frac{x}{|x|},\frac{I(x)}{|x|},\frac{J(x)}{|x|})$ is an orthonormal basis of $\bb{R}^3$. Then for $x,v,v_\ast\in\bb{R}^3$,
$\theta\in(0,\pi],~\varphi\in[0,2\pi)$, we set
\beq\label{para}
\left\{\begin{array}{l} \Gamma(x,\varphi):=(\cos\varphi) I(x)+(\sin\varphi) J(x),\\
v^\prime(v,v_\ast,\theta,\varphi):=v-\frac{1-\cos\theta}{2}(v-v_\ast)
+\frac{\sin\theta}{2}\Gamma(v-v_\ast,\varphi),\\
a(v,v_\ast,\theta,\varphi):=v^\prime(v,v_\ast,\theta,\varphi)-v,
\end{array}
\right.\eeq
then we write $\sigma\in\bb{S}^2$ as $\sigma=\frac{v-v_*}{|v-v_*|}\cos\theta + \frac{I(v-v_*)}{|v-v_*|}\sin\theta\cos\varphi + \frac{J(v-v_*)}{|v-v_*|}\sin\theta\sin\varphi$, and  observe at once that $\Gamma(x,\varphi)$ is orthogonal to $x$ and has the same norm as $x$,
from which it is easy to check that
\beq\label{num}
|a(v,v_*,\theta,\varphi)|=\sqrt{\frac{1-\cos\theta}{2}}|v-v_*|.
\eeq
Let us now give the definition  of weak solutions to \eqref{Bol}.
\bdefi\label{dfw}
Assume \eqref{con} is true for some $\nu\in(0,1), \gamma\in(-1,0)$ with $\gamma+\nu>0$.
A measurable family of probability measures $(f_t)_{t\geq 0}$
is called a \textit{weak solution} to \eqref{Bol} if it satisfies the following two conditions:
\begin{itemize}
  \item For all $t\geq 0$,
  \beq\label{conservation}
\int_{\bb{R}^3}vf_t(dv)=\int_{\bb{R}^3}vf_0(dv)~~\text{and}~~
\int_{\bb{R}^3}|v|^2f_t(dv)=\int_{\bb{R}^3}|v|^2f_0(dv)<\infty.
\eeq
  \item For any bounded globally Lipschitz function $\phi\in Lip(\bb{R}^3)$, any $t\in [0,T]$,
  \beq\label{weak}
  \int_{\bb{R}^3}\phi(v)f_t(dv)=\int_{\bb{R}^3}\phi(v)f_0(dv)+
  \int_0^t\int_{\bb{R}^3}\int_{\bb{R}^3}\ca{A}\phi(v,v_\ast)f_s(dv_\ast)f_s(dv)ds
  \eeq
  where
  \[\ca{A}\phi(v,v_\ast)=|v-v_*|^\gamma\int_{0}^{\pi/2}\beta(\theta) d\theta\int_0^{2\pi}
  [\phi(v + a(v,v_*,\theta,\varphi))-\phi(v)]d\varphi.\]
\end{itemize}
\edefi
We observe that $|\ca{A}\phi(v,v_*)|\le C_\phi|v-v_*|^{1+\gamma}\le C_\phi(1+|v-v_*|^2)$ from   $|a(v,v_*,\theta,\varphi)|\le C\theta|v-v_*|$ and $\int_0^{\pi/2}\theta\beta(\theta)d\theta<\infty$,   \eqref{weak} is thus well-defined.

\vskip1mm

Let us now recall the well-posedness result of
\eqref{Bol} in \cite[Corollary 2.4]{MR2511651} (more general existence results can be found in
\cite{MR1650006}).

\bthm\label{well-posedness} Assume \eqref{con}  for some $\gamma\in(-1,0)$,
 $\nu\in(0,1)$ with $\gamma+\nu>0$.  Let $q\ge2$ such that $q>\gamma^2/(\gamma+\nu)$.  Let
$f_0\in\ca{P}_q(\bb{R}^3)$ 
with $\int_{\bb{R}^3}f_0(v)|\log{f_0(v)}| dv< \infty$ and  let
$p\in(3/(3+\gamma),p_0(\gamma,\nu,q))$, where
\begin{equation}\label{pz}
p_0(\gamma,\nu,q)=\frac{q-\gamma}{q(3-\nu)/3-\gamma} \in (3/(3+\gamma),3/(3-\nu)).
\end{equation}
Then
\eqref{Bol} has a unique weak solution
$f\in L^\infty\big([0,\infty), \ca{P}_2({\bb{R}^3})\big)\cap L^1_{loc}\big([0,\infty), L^{p}({\bb{R}^3})\big)$.
\ethm

The explicit value of $p_0(\gamma,\nu,q)$ are not properly stated in \cite[Corollary 2.4]{MR2511651}.
However, following its proof (see the end of Step 3), we see that
$f\in L^1_{loc}\big([0,\infty), L^{p}({\bb{R}^3})\big)$ as soon as $1<p<3/(3-\nu)$ and
$-\gamma(p-1)/(1-p(3-\nu)/3)<q$. This precisely rewrites as
$p\in (1,p_0(\gamma,\nu,q))$.

\subsection{The particle system}
Let us now  recall  the Nanbu particle system introduced by \cite{nanbu1983interrelations}. It is
the $(\bb{R}^3)^N$-valued Markov process with infinitesimal generator $\ca{L}_N$ defined as
follows: for any bounded Lipschitz  function $\phi: (\bb{R}^3)^N\mapsto \bb{R}$ and
${\bf v}=(v_1,...,v_N)\in(\bb{R}^3)^N$,
$$
\ca{L}_{N}\phi({\bf v})=\frac{1}{N}\sum_{i\ne j}\int_{\bb{S}^2}[\phi({\bf v}+(v^\prime(v_i,v_j,\sigma)-v_i){\bf e}_i)-\phi({\bf v})]B(|v_i-v_j|, \theta)d\sigma,
$$
where $v{\bf e}_i=(0,...,v,...,0)\in(\bb{R}^3)^N$ with $v$ at the $i$-th place for $v\in\bb{R}^3$.

\vskip1mm

In other words, the system contains $N$ particles with velocities ${\bf v}=(v_1,...,v_N)$.
Each pair of particles (with velocities $(v_i,v_j)$), interact, for each $\sigma\in\bb{S}^2$,
at rate $B(|v_i-v_j|, \theta)/N$. Then one changes the velocity $v_i$ to
$v^\prime(v_i,v_j,\sigma)$ given by \eqref{Bol1} but $v_j$ remains unchanged.
That is, only one particle is changed at each collision.

\vskip1mm

The fact that  $\int_0^{\pi}\beta(\theta) d\theta=\infty$ (i.e. $\beta$ is non cutoff) means that
there are infinitely many jumps with a very small deviation angle.  It is thus impossible to
simulate it directly. For this reason, we will study a truncated version of Nanbu's particle system
applying a cutoff procedure as \cite{MR3456347}, who were studying the Nanbu system for
\emph{hard potentials} and \emph{Maxwell molecules}, and \cite{Cortez:2015aa}, who were dealing with
the Kac system for \emph{Maxwell molecules}. Our  particle system with cutoff
corresponds to the  generator $\ca{L}_{N,K}$ defined,
for any bounded Lipschitz  function $\phi: (\bb{R}^3)^N\mapsto \bb{R}$
and ${\bf v}=(v_1,...,v_N)\in(\bb{R}^3)^N$, by
\begin{align}\label{generator-cut}
\ca{L}_{N,K}\phi({\bf v})
&=\frac{1}{N}\sum_{i\ne j}\int_{\bb{S}^2}[\phi({\bf v}+(v^\prime(v_i,v_j,\sigma)-v_i){\bf e}_i)-\phi({\bf v})]
B(|v_i-v_j|, \theta)\bbd{1}_{\{\theta\ge G(K/|v_i-v_j|^\gamma)\}}d\sigma,
\end{align}
with $G$ defined by \eqref{nota}.

\vskip1mm

The generator $\ca{L}_{N,K}$ uniquely defines a strong Markov process with values
in $({\bb{R}^3})^N$. This comes from the fact that the corresponding jump rate
is finite and constant: for any configuration ${\bf v}=(v_1,...,v_N)\in(\bb{R}^3)^N$, it holds
that $N^{-1}\sum_{i\ne j}\int_{\bb{S}^2}B(|v_i-v_j|, \theta)\bbd{1}_{\{\theta\ge G(K/|v_i-v_j|^\gamma)\}}d\sigma
=2\pi (N-1)K$.
Indeed, for any $z\in [0,\infty)$, we have
$\int_{\bb{S}^2}B(x, \theta)\bbd{1}_{\{\theta\ge G(K/x^\gamma)\}}d\sigma=2\pi K$, which is easily checked
recalling that $B(x,\theta)=x^\gamma \beta(\theta)$ and the definition of $G$.

\subsection{Main results}
Now, we  give our uniqueness result for the Boltzmann equation.
\bthm\label{thm-uniqueness}
Assume \eqref{con} for some $\gamma\in(-1,0)$,  $\nu\in(0,1)$ satisfying  $\gamma+\nu>0$. Let $q\ge 2$
such that  $q>\gamma^2/(\gamma+\nu)$.  Assume that $f_0\in\ca{P}_q({\bb{R}^3})$ 
with a finite entropy,
i.e. $\int_{\bb{R}^3} f_0(v) |\log f_0(v)| dv<\infty$.  Let  $p\in (3/(3+\gamma),p_0(\gamma,\nu,q))$,
recall \eqref{pz},
and $(f_t)_{t\geq 0}\in L^\infty\big([0,\infty), \ca{P}_2({\bb{R}^3})\big)\cap L^1_{loc}\big([0,\infty), L^{p}({\bb{R}^3})\big)$ be the unique weak solution to \eqref{Bol} given by  Theorem \emph{\ref{well-posedness}}. Then for any other weak solution $(\tilde{f}_t)_{t\geq 0}\in L^\infty\big([0,\infty), \ca{P}_2({\bb{R}^3})\big)$ to  \eqref{Bol}, we have, for any $t\geq 0$,
\[
\ca{W}_2^2(f_t,\tilde{f}_t)\le \ca{W}_2^2(f_0,\tilde{f}_0)\exp{\Big(C_{\gamma,p}\int_0^t(1+\|f_s\|_{L^p})ds}\Big).\]
In particular,  we have uniqueness for \eqref{Bol} when starting from $f_0$ in the space of all weak solutions in the sense
of Definition \emph{\ref{dfw}}.
\ethm

The novelty of Theorem \ref{thm-uniqueness} is that no regularity at all is assumed
concerning $\tilde f$. In particular, we have uniqueness among all weak solutions,
while in \cite{MR2511651}, uniqueness is proved only in the class of weak solutions
lying in $L^\infty\big([0,\infty), \ca{P}_2({\bb{R}^3})\big)\cap L^1_{loc}\big([0,\infty),
L^{p}({\bb{R}^3})\big)$ for some $p>3/(3+\gamma)$.

\vskip1mm

Next, we  write  the following conclusion concerning  the particle system.
\bthm\label{main-result}
Assume \eqref{con} for some $\gamma\in(-1,0)$, $\nu\in(0,1)$ with $\gamma+\nu>0$.
Let $q>6$ such that  $q>\gamma^2/(\gamma+\nu)$ and let $f_0\in\ca{P}_q({\bb{R}^3})$
with a finite entropy.
Let $(f_t)_{t\ge0}$ be the  unique weak solution to  \eqref{Bol} given by Theorem \emph{\ref{well-posedness}}.
For each $N\ge 1$, $K\in[1,\infty)$, let  $(V_t^{i})_{i=1,...,N}$ be the Markov process with generator
$\ca{L}_{N,K}$ (see \eqref{generator-cut}) starting from an i.i.d. family $(V_0^{i})_{i=1,...,N}$ of
$f_0$-distributed random variables. We denote
the associated empirical measure by $\mu_{t}^{N,K}=N^{-1}\sum_{i=1}^N \delta_{V_t^{i}}$. Then for all $T>0$,
\[
\sup_{[0, T]}\bb{E}[\ca{W}_2^2(\mu_{t}^{N,K}, f_t)]\le C_{T,q}\Big(N^{-(1-6/q)(2+2\gamma)/3}+K^{1-2/\nu}+N^{-1/2}\Big).
\]
\ethm
We thus obtain a quantitive rate of chaos for the Nanbu's system with a singular interaction.
To our knowledge, this is the first result in this direction.
However, there is no doubt this rate is not the hoped optimal rate $N^{-1/2}$ like in the
hard potential case \cite{MR3456347}.

\subsection{Known results, strategies and main difficulties}

Let us give a non-exhaustive overview of the  known results on  the  well-posedness  of  \eqref{Bol}  for different potentials.  First,  the global existence of weak solution for the Boltzmann equation concerning all potentials was concluded by Villani in \cite{MR1650006}, with rather few assumptions on the initial  data (finite energy and entropy),
using some compactness methods.
However, the uniqueness results  are less well-understood. For  \emph{hard potentials} ($\gamma\in(0,1)$) with \emph{angular cutoff} ($\int_0^\pi \beta(\theta) d\theta<\infty$), there are some optimal results obtained by Mischler-Wennberg \cite{MR1697562}, where they  gave the existence of a unique weak $L^1$ solution to \eqref{Bol} with the minimal assumption that $\int_{\bb{R}^3}(1+|v|^2)f_0(v)dv<\infty$.
This was extended to weak measure solutions by Lu-Mouhot \cite{MR2871802}.
For the difficult case  \emph{without} angular cutoff, the first uniqueness result was obtained by Tanaka \cite{MR512334} concerning \emph{Maxwell molecules} ($\gamma=0$).
See also Toscani-Villani \cite{MR1675367}, who proved uniqueness
for  Maxwell molecules imposing that $\int_0^\pi \theta \beta(\theta) d\theta<\infty$
and that $\int_{\bb{R}^3}(1+|v|^2)f_0(dv)<\infty$.
Subsequently, Desvillettes-Mouhot \cite{MR2525118} (relying on a weighted $W_1^1$
space)  and Fournier-Mouhot \cite{MR2511651} (using the Wasserstein distance ${\cal W}_1$) successively  gave the uniqueness and stability  for both  \emph{hard potentials} ($\gamma\in(0,1]$) and \emph{moderately soft potentials} ($\gamma \in (-1,0)$ and $\nu \in (0,1)$) under different assumptions on initial data.  For \emph{moderately soft potentials}, the result in \cite{MR2511651} is much better since they use less assumptions on the initial condition than \cite{MR2525118}. Finally, let us mention another work \cite{MR2398952}, where Fournier-Gu\'{e}rin proved  a  local (in time) uniqueness result  with $f_0\in L^p(\bb{R}^3)$ for some $p>3/(3+\gamma)$ for the \emph{very soft potentials} ($\gamma \in (-3,0)$ and $\nu \in (0,2)$).

\vskip1mm

In this paper (Theorem \ref{thm-uniqueness}), we obtain a better uniqueness result in the case of a collision kernel without angular cutoff  when  $\gamma\in(-1,0)$ and $\nu\in(0,1-\gamma)$, that is, the uniqueness holds in the class of all measure solutions in $L^\infty\big([0,\infty), \ca{P}_2({\bb{R}^3})\big)$. This is very important when studying particle
systems. For example, a convergence result without rate would be almost immediate from our uniqueness: the tightness of the empirical measure of the particle system is not very difficult, as well as the fact that any limit point is a weak solution to \eqref{Bol}.
Since such a weak solution is unique by Theorem \ref{thm-uniqueness}, the convergence follows. Such a conclusion would be very difficult to obtain when using the uniqueness
proved in \cite{MR2511651}, because one would need to check that any limit point
of the empirical measure belongs to $L^1_{loc}([0,\infty,L^p(\bb{R}^3))$
for some $p>3/(3+\gamma)$,
which seems very difficult.

\vskip1mm

In order to extend the uniqueness result for all measure solutions, extra difficulty
is inevitable and the methods of \cite{MR2398952,MR2511651} will \emph{not} work.
However,  Fournier-Hauray  \cite{Fournier:2015aa} provide some ideas to overcome this,
in the simpler case of the Laudau equation for moderately soft potentials.
Here we follow these ideas, which rely on coupling methods.
Consider two weak solutions $f$ and $\tilde f$ in
$L^\infty\big([0,\infty), \ca{P}_2({\bb{R}^3})\big)$ to \eqref{Bol}, with possibly two
different initial conditions and assume that $f$ is {\it strong}, so that it belongs to
$L^1_{loc}\big([0,\infty), L^{p}({\bb{R}^3})\big)$.
First, we associate to the weak solution $\tilde f$
a {\it weak} solution $(X_t)_{t\geq 0}$ to some Poisson-driven SDE.
This uses a smoothing procedure as in
\cite{MR2375067,Fournier:2015aa}, but the situation is consequently
more complicated because we deal with jump processes.
Next, we try to associate to the strong solution $\tilde f$
a strong solution $(W_t)_{t\geq 0}$ to another SDE (driven by the same Poisson measure),
as \cite{Fournier:2015aa} did. But we did not manage to do this properly and we had
to use a  truncation procedure which though complicates our computation.
Then, roughly, we estimate ${\cal W}_2^2(f_t,\tilde f_t)$ by computing
$\E[|X_t-W_t|^2]$ as precisely as possible.

\vskip1mm 

The terminology \emph{propagation of chaos}, which is equivalent to the convergence of the empirical measure of a  particle system  to the solution to a nonlinear equation, was first formulated by  Kac \cite{MR0084985}. He was studying the convergence of a {\it toy} particle system as a step to the rigorous derivation of the Boltzmann equation.
Afterwards,  McKean \cite{MR0224348} and Gr\"{u}nbaum \cite{MR0334788} extended Kac's ideas to study the chaos property for different models with bounded collision kernels.  Sznitman \cite{MR753814} then showed the chaos property (without rate) for the \emph{hard spheres} ($\gamma=1$ and $\nu=0$). Following  Tanaka's  probabilistic interpretation for the Boltzmann equation with \emph{Maxwell molecules}, Graham-M\'{e}l\'{e}ard \cite{MR1428502}  were the first to give a rate of chaos for \eqref{Bol}, concerning both Kac and Nanbu models, for Maxwell molecules with cutoff ($\gamma=0$ and $\int_0^\pi \beta(\theta)d\theta<\infty$), using the  total variation distance. Recently,  some important progresses have
been made. First,  Mischler-Mouhot \cite{MR3069113} obtained a uniform (in time) rate of convergence of Kac's  particle system of order $N^{-\e}$ (for Maxwell molecules without cutoff) and $(\log N)^{-\e}$ (for hard spheres, i.e. $\gamma=1$ and $\nu=0$), with some
small $\e>0$.
This result, entirely relying on analytic methods, is noticeable, although the rates
are clearly not sharp.
Then, Fournier-Mischler \cite{MR3456347} proved the propagation of chaos at rate
$N^{-1/4}$ for the Nanbu system and for hard potentials without cutoff
($\gamma \in [0,1]$ and $\nu \in (0,1)$).
Finally,  as mentioned in Section 1.5,  Cortez-Fontbona \cite{Cortez:2015aa} used
two coupling techniques for Kac's binary interaction system and  obtained
a uniform in time estimate for the Boltzmann equation with \emph{Maxwell molecules}
($\gamma=0$) under some suitable moments assumptions on the initial datum.
Let us mention that the time-uniformity uses the recent nice results of
Rousset \cite{Rousset:2014aa}.

\vskip1mm

In this paper (Theorem \ref{main-result}),  we obtain, to our knowledge, the first chaos result (with rate)
for soft potentials (which are, of course, more difficult), but it is a bit
unsatisfying: (1) we cannot study  Kac's system (which is physically more reasonnable
than Nanbu's system) because it is not readily to exhibit a suitable coupling;
(2) our consideration is merely for $\gamma\in(-1,0)$, since some basic estimates
in Section 2 do not hold any more if $\gamma \le -1$;
(3) our rate is not sharp.
However, since the interaction is  singular, it seems hopeless to get a perfect result.

\vskip1mm

In terms of the propagation of chaos with a singular interaction,  there are only very few results.  Hauray-Jabin \cite{MR3377068} considered a deterministic system  of  particles interacting  through a force of the type $1/|x|^\alpha$ with $\alpha< 1$,
in dimension $d \ge 3$, and proved the mean field limit and the propagation of chaos
to the Vlasov equation. Also, Fournier-Hauray-Mischler \cite{MR3254330} proved the convergence of the vortex model to the 2D Navier-Stokes equation with a singular Biot-Savart kernel using some entropy dissipation technique.  Following the method of \cite{MR3254330}, Godinho-Qui\~{n}inao \cite{MR3365970} proved the propagation of chaos of some particle system to the 2D subcritical Keller-Segel equation.
Recently, Fournier-Hauray \cite{Fournier:2015aa} proved  propagation of chaos for the Landau equation with a singular interaction ($\gamma\in(-2,0)$). Actually, they gave a quantitative rate of chaos when $\gamma\in(-1,0)$, while the convergence without rate was checked
when $\gamma\in(-2,0)$ by the entropy dissipation technique.

\vskip1mm

Roughly speaking,  to prove  our propagation of chaos result,
we consider an approximate version of our stability principle,
with a discrete $L^p$ norm as in \cite{Fournier:2015aa}. Here, we list the main
difficulties: The trajectory of a typical particle related to the Boltzmann
equation is a jump process so that all the continuity arguments used in
\cite{Fournier:2015aa} have to be changed. In particular, a detailed
study of small and large jumps is required.
Also, the solution to the Landau equation lies in $L^1_{loc}\big([0,\infty),L^2(\bb{R}^3)\big)$,
while the one of the Boltzmann equation  lies in $L^1_{loc}\big([0,\infty),L^p(\bb{R}^3)\big)$
for some $p$ smaller than $2$. This causes a few difficulties in Section \ref{ttt},
because working in $L^p$ is slightly more complicated.

\subsection{Arrangement of the paper and final notations}
In Section 2,  we give some basic estimates. In Section 3, we establish the strong/weak stability principle for \eqref{Bol}. In Section 4, we construct the suitable coupling.
In Section 5, we bound the $L^p$ norm of an empirical measure in terms of  $L^p$ norm of  the weak solution.
Finally, in Section 6, we prove the convergence of the particle system.

\vskip 1mm

In the  sequel, $C$ stands for a positive constant whose value may change from line to line.
When necessary, we will indicate in subscript the parameters it depends on.

\vskip1mm

In the whole paper, we consider two probability spaces by  Tanaka's idea for the probabilistic interpretation of the Boltzmann equation in Maxwell molecules case:  the first space is the abstract space $(\Omega, \ca{F}, \bb{P})$ and the second  is $([0,1],\ca{B}([0,1]),d\alpha)$.  A stochastic process defined on the latter space is called an $\alpha$-processes and we denote the expectation on $[0,1]$ by $\bb{E}_\alpha$ and the laws by $\ca{L}_\alpha$.

\section{Preliminaries}

Above all, let us recall that  for $\gamma\in(-1,0)$, $p>3/(3+\gamma)$ and
$f\in\ca{P}(\bb{R}^3)\cap L^p(\bb{R}^3)$, it holds that
\begin{align}\label{norm-inequality}
 \sup_{v\in\bb{R}^3}\int_{\bb{R}^3}|v-v_*|^\gamma f(dv_*)
 &\le \sup_{v\in\bb{R}^3}\int_{|v-v_*|\le1}|v-v_*|^\gamma f(dv_*) + \sup_{v\in\bb{R}^3}\int_{|v-v_*|\ge1}|v-v_*|^\gamma f(dv_*)\nonumber\\
 &\le 1+ C_{\gamma,p}\|f\|_{L^{p}(\bb{R}^3)},
\end{align}
where $C_{\gamma,p}=\sup_{v\in\bb{R}^3}[\int_{|v-v_*|\le1}|v-v_*|^{p\gamma/(p-1)}dv_*]^{(p-1)/p}=[\int_{|v_*|\le1}|v_*|^{p\gamma/(p-1)}dv_*]^{(p-1)/p}<\infty$, since $p>3/(3+\gamma)$ by assumption.

\vskip1mm

Let us  now classically rewrite  the collision operator by making disappear the velocity-dependence $|v-v_*|^{\gamma}$ in the \emph{rate} using a substitution.
\blem\label{llxx}
We assume \eqref{con} and recall \eqref{nota} and \eqref{para}. For $z\in [0,\infty)$, $\varphi\in[0,2\pi)$, $v,v_*\in\bb{R}^3$ and $K\in[1,\infty)$, we define
\beq\label{rewrite}
c(v,v_*,z,\varphi):=a[v,v_*,G(z/|v-v_*|^{\gamma}),\varphi]~~~\text{and} ~~~c_K(v,v_*,z,\varphi):=c(v,v_*,z,\varphi)\bbd{1}_{\{z\le K\}}.
\eeq
For any $\phi\in Lip(\bb{R}^3)$, any $v,v_*\in\bb{R}$,
\beq\label{opera}
\ca{A}\phi(v,v_*)=\int_0^\infty dz\int_0^{2\pi}d\varphi [\phi(v+c(v,v_*,z,\varphi))-\phi(v)].
\eeq
For any $N\ge 1$, $K\in[1,\infty)$, ${\bf v}=(v_1,...,v_N)\in(\bb{R}^3)^N$, any bounded measurable $\phi:
(\bb{R}^3)^N\mapsto\bb{R}$,
\beq\label{rewpartopera}
\ca{L}_{N,K}\phi({\bf v})=\frac{1}{N}\sum_{i\ne j}\int_0^\infty
dz\int_0^{2\pi}d\varphi[\phi({\bf v}+c_K(v_i,v_j,z,\varphi){\bf e}_i)-\phi({\bf v})].
\eeq
\elem
This lemma is stated in \cite[Lemma 2.2]{MR3456347} when $\gamma\in[0,1]$, but the proof does not use this fact:
it actually holds true for any $\gamma\in\bb{R}$.
Next, let us recall Lemma 2.3 in \cite{MR3456347}  which is an accurate version of Tanaka's trick in \cite{MR512334}. Here, we adopt  the notation \eqref{para}.
\blem\label{Tanaka}
There exists some 
measurable function $\varphi_0: \bb{R}^3\times\bb{R}^3\mapsto[0,2\pi)$
such that for all $X, Y\in\bb{R}^3$, all $\varphi\in[0,2\pi)$,
\begin{align*}
|\Gamma(X,\varphi) - \Gamma(Y,\varphi+\varphi_0(X,Y)|\le|X-Y|.
\end{align*}
\elem
The rest of the section  is  an adaption of \cite[Section 3]{MR3456347},
which assumes that $\gamma\in[0,1]$, to the case where $\gamma\in(-1,0)$.
When compared with  \cite{MR2398952}, what is new is that in the inequalities \eqref{ee2} and \eqref{ee3} below,
only $|v-v_*|^\gamma$ appears (while in \cite{MR2398952}, there is $|v-v_*|^\gamma+|\tilde v-\tilde v_*|^\gamma$).
This is very useful to get a strong/weak stability estimate:
we will be able to use the regularity of only one of the two solutions to be compared.
Let us mention that it seems impossible to extend our ideas to the more singular case where $\gamma \leq -1$.

\blem\label{estimategeneral}
There is a constant $C$ such that for any $v, v_*,\tilde{v}, \tilde{v}_*\in\bb{R}^3$, any $K\geq 1$,
\begin{align}
\int_0^{\infty}\int_0^{2\pi}|c(v,v_*,z,\varphi)&-c(\tilde{v}, \tilde{v}_*,z,\varphi+\varphi_0(v-v_*, \tilde{v}-\tilde{v}_*))|^2d\varphi dz\label{ee2}\\
\leq& C(|v-\tilde{v}| ^2 + |v_*-\tilde{v}_*|^2) |v-v_*|^\gamma,\notag\\
\int_0^{\infty}\int_0^{2\pi}\Big(|v+c(v,v_*,z,\varphi)-&\tilde{v}-c_K(\tilde{v}, \tilde{v}_*,z,\varphi+\varphi_0(v-v_*, \tilde{v}-\tilde{v}_*))|^2-|v-\tilde{v}|^2 \Big)d\varphi dz\label{ee3}\\
\leq& C (|v-\tilde{v}|^2+|v_*-\tilde{v}_*|^2 )|v-v_*|^\gamma + C |v-v_*|^{2+2\gamma/\nu}K^{1-2/\nu} .\notag\\
\int_0^\infty \int_0^{2\pi}|c_K(v,v_*,z,\varphi)|^2 d\varphi dz \leq &C|v-v_*|^{\gamma+2}, \quad
 \int_0^\infty \left | \int_0^{2\pi} c_K(v,v_*,z,\varphi) d\varphi\right | dz \le C|v-v_*|^{\gamma+1}\label{ee4}\\
 \int_0^\infty \int_0^{2\pi}|c(v,v_*,z,\varphi)|^2 d\varphi dz \leq &C|v-v_*|^{\gamma+2}, \quad
 \int_0^\infty \left | \int_0^{2\pi} c(v,v_*,z,\varphi) d\varphi\right | dz \le C|v-v_*|^{\gamma+1}\label{ee5}
\end{align}
\elem

\bpf
For $x>0$, we set $\Phi_K(x)= \pi\int_0^{K}(1-\cos G(z/x^\gamma))dz$  and
$\Psi_K(x)= \pi\int_{K}^\infty(1-\cos G(z/x^\gamma))dz$.
We introduce the shortened notation
$x=|v-v_*|$, $\tilde x=|\tilde v-\tilde v_*|$, $\varphi_0=\varphi_0(v-v_*, \tilde{v}-\tilde{v}_*)$,
$c=c(v,v_*,z,\varphi)$, $c_K=c_K(v,v_*,z,\varphi)=c\bbd{1}_{\{z\leq K\}}$,
$\tilde c=c(\tilde{v}, \tilde{v}_*,z,\varphi+\varphi_0)$ and
$\tilde c_K=c_K(\tilde{v}, \tilde{v}_*,z,\varphi+\varphi_0)=\tilde c\bbd{1}_{\{z\le K\}}$.

\vskip1mm

{\it Step 1.} We first verify that $\Phi_K(x) \leq C x^\gamma$ and that $|\Phi_K(x)-\Phi_K(\tilde x)|\leq
C |x^\gamma-\tilde x^\gamma|$. First, we immediately see that $\Phi_K(x)\leq \pi \int_0^\infty G^2(z/x^\gamma)dz
=x^\gamma \pi \int_0^\infty G^2(z)dz$ which implies the first point (recall \eqref{con1}).
To check the second point, it suffices to verify that $F_K(x)=\int_0^K(1-\cos G(z/x))dz$
has a bounded derivative (uniformly in $K\geq 1$). But
we have $F_K(x)=x \int_0^{K/x}(1-\cos G(z))dz$ so that
\begin{align*}
 |F_K^\prime(x)|
 \le  \int_0^\infty (1-\cos G(z)) dz+ x (K/x^2)(1-\cos G(K/x))
 \le  C+ (K/x) G^2(K/x),
 \end{align*}
which is uniformly bounded by \eqref{con1}.

\vskip1mm

{\it Step 2.} Proceeding as in the proof of  \cite[Lemma 3.1]{MR3456347}, we see that
$\int_0^\infty \int_0^{2\pi}|c_K|^2 d\varphi dz=x^2\Phi_K(x)$, which is bounded by $Cx^{\gamma+2}$
by Step 1. Also, recalling \eqref{para} and \eqref{rewrite}, using that $\int_0^{2\pi} \Gamma(X,\varphi)d\varphi=0$,
we see that $\int_0^{2\pi} c_K d\varphi= -\pi (v-v_*) (1-\cos G(z/x^\gamma))$,
whence $\int_0^\infty |\int_0^{2\pi} c_K d\varphi |dz=x \Phi_K(x) \leq C x^{\gamma+1}$ by Step 1.
All this proves \eqref{ee4}, from which \eqref{ee5} follows by letting $K$ increase to infinity.

\vskip1mm

{\it Step 3.}
Let us denote by $I_K=\int_0^{K}\int_0^{2\pi} |c - \tilde c|^2 d\varphi dz $, by
$J_K=\int_0^K\int_0^{2\pi} (|v+c-\tilde v - \tilde c|^2-|v-\tilde v|^2) d\varphi dz$
and by $L_K=\int_K^\infty\int_0^{2\pi} (|v+c-\tilde v|^2-|v-\tilde v|^2) d\varphi dz$.
Proceeding exactly as in the proof of \cite[Lemma 3.1]{MR3456347}, we see that
$J_K \leq A_1^K+A_2^K$ and $L_K \leq A^K_3$, where
\begin{align*}
&A_1^K=2 x \tilde x \int_0^{K}\Big(G(z/x^\gamma) - G(z/\tilde{x}^\gamma)\Big)^2 dz,\\
&A_2^K=\big[|v-\tilde{v}|+|v_*- \tilde{v}_*|\big]|(v-v_*)\Phi_K(x)-(\tilde{v} - \tilde{v}_*)\Phi_K(
\tilde{x})|,\\
&A_3^K=(x^2+2|v-\tilde{v}| x)\Psi_K(x).
\end{align*}

Also, $I_K=J_K -2 (v-\tilde v) \cdot \int_0^{K}\int_0^{2\pi} (c - \tilde c) d\varphi dz $ and, as seen in
the proof of \cite[Lemma 3.1]{MR3456347}, $\int_0^K\int_0^{2\pi} c d\varphi dz = -(v-v_*)\Phi_K(x)$, so that
$I_K \leq J_K+A_4^K$ with
$$
A_4^K = 2 |v-\tilde v| |(v-v_*)\Phi_K(x)-(\tilde v-\tilde v_*)\Phi_K(\tilde x) |.
$$

First, we immediately deduce from  \eqref{ineq} that
\begin{align*}
A_1^K&\le 2c_4 x \tilde{x} \frac{(x^\gamma - \tilde{x}^\gamma)^2}{x^\gamma + \tilde{x}^\gamma}
\le 2 c_4 (x-\tilde x)^2 \min{(x^\gamma, \tilde{x}^\gamma)}
\leq C (|v-\tilde{v}| ^2 + |v_*-\tilde{v}_*|^2)  |v-v_*|^\gamma.
\end{align*}
For the second inequality, we used that $|x^{\gamma} - \tilde x^{\gamma}|\le|x^{-1}-\tilde x^{-1}|(x\land \tilde x)^{1+\gamma}$
(because $\gamma \in (-1,0)$) so that
\begin{align*}
x\tilde x\frac{|x^\gamma - \tilde x^\gamma|^2}{x^\gamma+\tilde x^\gamma}
 \le (x\tilde x)^{1+|\gamma|}\frac{|x^{-1}-\tilde x^{-1}|^2 (x\land \tilde x)^{2\gamma+2}}{x^{|\gamma|}+\tilde x^{|\gamma|}}
\le (x\tilde x)^{|\gamma|-1}\frac{|x-\tilde x|^2 (x\tilde x)^{1+\gamma}}{x^{|\gamma|}+\tilde x^{|\gamma|}} =
\frac{|x-\tilde x|^2} {x^{|\gamma|}+\tilde x^{|\gamma|}},
\end{align*}
which is indeed bounded by $(x-\tilde x)^2 \min{(x^\gamma, \tilde{x}^\gamma)}$.

\vskip1mm

We now verify that
$A_2^K \leq C\big(|v-\tilde{v}|^2 + |v_*- \tilde{v}_*|^2\big)|v-v_*|^\gamma$.
By Step 1, for any $X,Y\in\bb{R}^3$,
\[\left|X\Phi_K(|X|)- Y\Phi_K(|Y|)\right|\le |Y||\Phi_K(|X|)-\Phi_K(|Y|)|+|X-Y|\Phi_K(|X|)\le
C |Y|\Big ||X|^\gamma-|Y|^\gamma\Big|+C|X-Y| |X|^\gamma.\]
Since again $|x^{\gamma} - \tilde x^{\gamma}|\le|x^{-1}-\tilde x^{-1}|(x\land \tilde x)^{1+\gamma}$,
we conclude that
 $\left|X\Phi_K(|X|)- Y\Phi_K(|Y|)\right| \leq C |X-Y| |X|^\gamma$, whence
\begin{align*}
A_2^K &\le C\, \big[|v-\tilde{v}|+|v_*- \tilde{v}_*|\big]|(v-v_*)-(\tilde{v} - \tilde{v}_*)|\min\{x^\gamma, \tilde x^\gamma \}
\end{align*}
as desired.

\vskip1mm

We next observe that $A_4^K \leq 2A_2^K$.

\vskip1mm

Finally, we see that $\Psi_K(x)\le C\int_K^\infty G^2(z/x^\gamma)dz\le
C\int_K^\infty(z/x^\gamma)^{-2/\nu}dz=Cx^{2\gamma/\nu}K^{1-2/\nu}$  and that $\Psi_K(x)\le C\int_0^\infty G^2(z/x^\gamma)dz\le C\int_0^\infty (1+z/x^\gamma)^{-2/\nu}dz =C x^\gamma$ according to  \eqref{con1} , which imply  $\Psi_K(x)\le C\min\{x^\gamma, x^{2\gamma/\nu}K^{1-2/\nu}\}$. Hence,
\[
A_3^K=(x^2+2|v-\tilde v| x ) \Psi_K(x) \leq C |v-\tilde{v}|^2 |v-v_*|^\gamma
+ C |v-v_*|^{2+2\gamma/\nu}K^{1-2/\nu},
\]
because $2|v-\tilde v| x \le |v-\tilde v|^2 + x^2$ and $x^2\Psi_K(x)\le C x^{2+2\gamma/\nu}K^{1-2/\nu}$.
%and $2|v-\tilde{v}| x\Psi_K(x)\le |v-\tilde{v}|^2 x^\gamma+x^{2-\gamma}\Psi^2_K(x)$
%with $\Psi^2_K(x) \leq C x^{ \gamma + 2\gamma/\nu}K^{1-2/\nu}$.

\vskip1mm

The left hand side of \eqref{ee3} is nothing but
$J_K+L_K$, which is bounded by $A_1^K+A_2^K+A_3^K$: \eqref{ee3} is proved.
Finally, the left hand side of \eqref{ee2} equals $\lim_{K\to \infty} I_K$ and we know that
$I_K \leq A_1^K+A_2^K+A_4^K$, which is (uniformly in $K$) bounded by
$(|v-\tilde{v}|^2 + |v_*- \tilde{v}_*|^2) |v-v_*|^\gamma$ as desired.
\epf

\section{Stability}
In this section, our goal is  to prove Theorem \ref{thm-uniqueness}.
For this, we first need two important propositions.
The first one is the most delicate part in the whole proof, which consists in showing that we can associate a {\it weak}
solution to some SDE to {\it any} weak solution to \eqref{Bol}.
It extends Proposition B.1 in \cite{Fournier:2015aa} to our situation, see also
\cite[Theorem 2.6]{MR2375067} (both concerning diffusion-type PDEs and Brownian SDEs). We will prove it later.
\bprop\label{prop-PDE}
Assume \eqref{con} for some $\gamma\in(-1,0)$,  $\nu\in(0,1)$ with $\gamma+\nu>0$.
Consider any weak solution $(\tilde f_t)_{t\ge0}\in L^\infty\big([0,\infty), \ca{P}_2({\bb{R}^3})\big)$ to
\eqref{Bol}. Then there exists, on some probability space, a random variable $X_0$ with law $\tilde f_0$, independent of
a Poisson measure
$M(ds, d\alpha,dz, d\varphi)$ on $[0, \infty) \times [0,1]\times [0,\infty)\times[0,2\pi)$ with intensity
$ds d\alpha dz d\varphi$, a measurable family $(X^*_t)_{t\geq 0}$ of $\alpha$-random variables and a c\`adl\`ag
adapted process $(X_t)_{t\ge0}$ solving
\beq\label{SDE}
 X_t=X_0+\int_0^t\int_0^1\int_0^\infty\int_0^{2\pi} c \big(X_{s-}, X_s^*(\alpha), z,\varphi\big)
M(ds,d\alpha, dz,d\varphi)
\eeq
and such that for all $t\geq 0$, $\cL(X_t)=\cL_\alpha(X^*_t)=\tilde f_t$.
\eprop
The second one is the following proposition.

\bprop\label{prop-PDE2}
Assume \eqref{con} for some $\gamma\in(-1,0)$,  $\nu\in(0,1)$ with $\gamma+\nu>0$,
that $f_0\in \ca{P}_q(\bb{R}^3)$ for some $q\ge2$ such that $q>\gamma^2/(\gamma+\nu)$ and that
$f_0$ has a finite entropy. Fix $p\in(3/(3+\gamma),p_0(\gamma,\nu,q))$.
Let $(f_t)_{t\ge0} \in L^\infty\big([0,\infty), \ca{P}_2({\bb{R}^3})\big)
\cap L^1_{loc}\big([0,\infty), L^{p}({\bb{R}^3})\big)$ be the corresponding unique weak solution to  \eqref{Bol} given by Theorem \emph{\ref{well-posedness}}.
Consider also the Poisson measure $M$, the process $(X_t)_{t\geq 0}$ and the family $(X^*_t)_{t\geq 0}$ built in
Proposition \emph{\ref{prop-PDE}} (associated to another weak solution $(\tilde f_t)_{t\ge0}\in L^\infty\big([0,\infty),
\ca{P}_2({\bb{R}^3})\big)$. Let $W_0\sim f_0$ (independent of $M$) be such that
$\E[|W_0-X_0|^2]=\cW_2^2(f_0,\tilde f_0)$ and, for each $t\geq 0$, an $\alpha$-random variable $W^*_t$
such that $\cL_\alpha(W^*_t)=f_t$ and $\E_\alpha[|W_t^*-X_t^*|^2]=\cW_2^2(f_t,\tilde f_t)$.
Then for $K \ge 1$, the equation
\beq\label{Bprocess}
 W_t^K=W_0+\int_0^t\int_0^1\int_0^\infty\int_0^{2\pi}c_K(W_{s-}^K, W_s^*(\alpha),z,\varphi+\varphi_{s,\alpha,K})M(ds,d\alpha,dz,d\varphi),
\eeq
with $\varphi_{s,\alpha,K}=\varphi_0(X_{s-}-X_s^*(\alpha),W_{s-}^K-W_s^*(\alpha))$,
has a unique solution. Moreover, setting $f_t^K=\cL(W^K_t)$ for each $t\ge 0$, it holds that for all $T>0$,
\beq\label{prop-law-convergence}
\lim_{K\rightarrow\infty}\sup_{[0,T]}\ca{W}_2^2(f_t^K, f_t)=0.
\eeq

\eprop

\bpf
For any $K\geq 1$,
the Poisson measure involved in \eqref{Bprocess} is actually finite
(because $c_K=c\bbd{1}_{\{z\le K\}}$), so the existence and uniqueness for this equation is obvious.
It only remains to prove \eqref{prop-law-convergence}, which has already been done
in \cite[Lemma 4.2]{MR2398952}, where the formulation of the equation
is slightly different. But one easily checks that $(W^K_t)_{t\geq 0}$ is a (time-inhomogeneous)
Markov process with the same generator as the one defined by \cite[Eq. (4.1)]{MR2398952},
because for all bounded measurable function
$\phi:\bb{R}^3\mapsto \bb{R}$ and all $t\geq 0$, a.s.,
\begin{align*}
&\int_0^1\int_0^\infty\int_0^{2\pi}\Big[\phi(w+c_K(w, W_t^*(\alpha),z,\varphi+
\varphi_0(X_{t-}-X^*_t(\alpha),w-W^*_t(\alpha)) )-\phi(w)\Big]  d\varphi dz d\alpha\\
=& \int_0^1\int_0^\infty\int_0^{2\pi}\Big[\phi(w+c_K(w,v,z,\varphi))-\phi(w)\Big]
d\varphi dz f_t(dv)
\end{align*}
by the $2\pi$-periodicity of $c_K$ (in $\varphi$) and since  $\cL_\alpha(W^*_t)=f_t$.
\epf

Now, we use these coupled processes to conclude the 
\bpf[Proof of Theorem \ref{thm-uniqueness}]
We consider a weak solution $(\tilde f_t)_{t\geq 0}$ to \eqref{Bol}, with which we associate the objects
$M$, $(X_t)_{t\geq 0}$, $(X_t^*)_{t\geq 0}$ as in Proposition \ref{prop-PDE}. We then consider $f_0$
satisfying the assumptions of Theorem \ref{well-posedness} and the corresponding unique weak solution
$(f_t)_{t\geq 0}$ belonging to  $L^\infty\big([0,\infty), \ca{P}_2({\bb{R}^3})\big)
\cap L^1_{loc}\big([0,\infty), L^{p}({\bb{R}^3})\big)$ (with $p\in(3/(3+\gamma),p_0(\gamma,\nu,q))$)
and we consider $(W_t^K)_{t\geq 0}$, $(W_t^*)_{t\geq 0}$
built in Proposition \ref{prop-PDE2} for any $K\ge 1$. We know that
$\cW_2^2(f_0,\tilde f_0)=\E[|W_0-X_0|^2]$ and that
$\cW_2^2(f_t,\tilde f_t)=\E_\alpha[|W_t^*-X_t^*|^2]$
for all $t\geq 0$.  Using  that $W_t^K\sim f_t^K$ and $X_t \sim \tilde f_t$ for each $t\ge0$,
we deduce from  \eqref{prop-law-convergence} that for all $t\ge0$,
\begin{equation}\label{dfjt}
\ca{W}_2^2(f_t,\tilde{f}_t)\le \limsup_{K\to\infty}\bb{E}[|W_t^K-X_t|^2]=:J_t.
\end{equation}
 Next, we focus on the time interval $[0, T]$ for any fixed $T>0$, and  split the proof into several steps.

\vskip1mm

{\it Step 1.}
By the It\^{o} formula,  we know  that
\begin{align*}
\bb{E}[|W_t^K-X_t|^2]=\bb{E}[|W_0-X_0|^2]+\bb{E}\left[ \int_0^t \int_0^1 \Delta_s^K(\alpha) d\alpha ds \right],
\end{align*}
where
\begin{align*}
\Delta_s^K(\alpha):=&\int_0^\infty\int_0^{2\pi}\Big(|W_{s}^K-X_{s}+c_{K,W}(s) - c_X(s)|^2-|W_{s}^K-X_{s}|^2\Big) d\varphi dz
\end{align*}
with the shortened notation $c_{K,W}(s):=c_K \big(W_{s}^K,W_s^*(\alpha),z,\varphi+\varphi_{s,\alpha,K}\big)$
and $c_X(s):=c \big(X_{s},X_s^*(\alpha),z,\varphi\big).$
We then show that
\begin{align}
\Delta_s^K(\alpha) \le& C (|W_s^K-X_s| ^2 + |W_s^*(\alpha)-X_s^*(\alpha)|^2) |W_s^K - W_s^*(\alpha)|^\gamma
+ C |W_s^K - W_s^*(\alpha)|^{2+2\gamma/\nu} K^{1-2/\nu}, \label{rrr1}
\end{align}
and
\begin{align}
\Delta_s^K(\alpha) \leq & C|W_s^K-W_s^*(\alpha)|^{\gamma+2}+C|X_s-X_s^*(\alpha)|^{\gamma+2} \notag \\
&+C |W_{s}^K-X_{s}|\big(|W_s^K-W_s^*(\alpha)|^{\gamma+1}+|X_s-X_s^*(\alpha)|^{\gamma+1}\big).\label{rrr2}
\end{align}
First, Lemma \ref{estimategeneral} (inequality \eqref{ee3}) precisely tells us that \eqref{rrr1} holds true.
Next, we observe that
\begin{align*}
\Delta_s^K(\alpha)
\le  2\int_0^\infty\int_0^{2\pi} (|c_{K,W}(s)|^2+|c_X(s)|^2)  d\varphi dz+ 2|W_{s}^K-X_{s}|\Big|\int_0^\infty\int_0^{2\pi}(c_{K,W}(s) - c_X(s)) d\varphi dz \Big|.
\end{align*}
Hence,  using \eqref{ee4} and \eqref{ee5}, the proof of \eqref{rrr2} is concluded.

\vskip1mm 

{\it Step 2.} Set $\kappa(\gamma)=\min((\gamma+1)/|\gamma|,|\gamma|/2)>0$.  We verify that  there exists a constant $C(T, f_0, \tilde{f}_0, f)>0$ (depending on $T$, $m_2(f_0)$, $m_2(\tilde{f}_0)$, $\int_0^t \|f_s\|_{L^{p}}ds$),  such that for all $\ell \ge 1$ (and all $K\geq 1$),
$$
I_t^{i,\ell} \le C(T, f_0, \tilde{f}_0, f)\ell^{-\kappa(\gamma)}, \quad  i=1,2,3,4,
$$
where
\begin{align*}
&I_t^{1,\ell}:=\bb{E}\Big[\int_0^t \int_0^1|W_s^K-W_s^*(\alpha)|^{\gamma+2}\,\bbd{1}_{\{|W_s^K - W_s^*(\alpha)|^\gamma\ge \ell\}} d\alpha ds\Big],\\
&I_t^{2, \ell}:=\bb{E}\Big[\int_0^t \int_0^1|X_s-X_s^*(\alpha)|^{\gamma+2}\, \bbd{1}_{\{|W_s^K - W_s^*(\alpha)|^\gamma\ge \ell\}} d\alpha ds\Big],\\
&I_t^{3,\ell}:=\bb{E}\Big[\int_0^t \int_0^1|W_{s}^K-X_{s}||W_s^K-W_s^*(\alpha)|^{\gamma+1}\,\bbd{1}_{\{|W_s^K - W_s^*(\alpha)|^\gamma\ge \ell\}} d\alpha ds\Big],\\
&I_t^{4, \ell}:=\bb{E}\Big[\int_0^t \int_0^1|W_{s}^K-X_{s}||X_s-X_s^*(\alpha)|^{\gamma+1}\, \bbd{1}_{\{|W_s^K - W_s^*(\alpha)|^\gamma\ge \ell\}} d\alpha ds\Big].
\end{align*}

Since $\gamma\in(-1,0)$ and $\kappa(\gamma)\le (\gamma+2)/|\gamma|$, we have
\[
I_t^{1,\ell}\leq \ell^{-(\gamma+2)/|\gamma|}T\le \ell^{-\kappa(\gamma)}T.
\]
Similarly,
\begin{align*}
I_t^{3,\ell} &\le \ell^{-(\gamma+1)/|\gamma|} \int_0^t \bb{E}\Big[ |W_{s}^K-X_{s}|\Big] ds.
\end{align*}
Using \eqref{conservation} for $(f_t)_{t\ge0}$ and $(\tilde{f}_t)_{t\ge0}$, \eqref{prop-law-convergence},  and that $m_2(f_s^K)\le 2m_2(f_s)+2\ca{W}_2^2(f_s,f_s^K)$, we know that  $\bb{E}\Big[ |W_{s}^K-X_{s}|\Big]\le
C(1+m_2(f_s^K)+m_2(\tilde{f}_s))\le C(T, f_0, \tilde{f}_0)$. Hence,
\[I_t^{3,\ell}\le C(T, f_0, \tilde{f}_0)\ell^{-\kappa(\gamma)}. \]
Since $\gamma+2\in(1,2)$, it follows from the H\"{o}lder inequality that
\begin{align*}
I_t^{2, \ell}
&\le \bb{E}\left[\Big(\int_0^t \int_0^1|X_s-X_s^*(\alpha)|^{2}d\alpha ds\Big)^{\frac{\gamma+2}{2}}\Big( \int_0^t \int_0^1\bbd{1}_{\{|W_s^K - W_s^*(\alpha)|^\gamma\ge \ell\}} d\alpha ds\Big)^{\frac{|\gamma|}{2}}\right]\\
&\le C\bb{E}\left[\Big(\int_0^t (|X_s|^2+m_2(\tilde{f}_s)) ds\Big)^{\frac{\gamma+2}{2}}\Big( \int_0^t \int_0^1\frac{|W_s^K - W_s^*(\alpha)|^\gamma}{\ell} d\alpha ds\Big)^{\frac{|\gamma|}{2}}\right]
\end{align*}
Since $\cL_\alpha(W_s^*)=f_s$, we have
$\int_0^1|W_s^K - W_s^*(\alpha)|^\gamma d\alpha= \int_{\bb{R}^3} |W_s^K - v|^\gamma f_s(dv)\leq
1+ C_{\gamma,p}\|f_s\|_{L^{p}}$ by \eqref{norm-inequality},
so that
\begin{align*}
I_t^{2, \ell}
\le&  \ell^{\gamma/2}\Big(1+ \int_0^t \big(\bb{E}[|X_s|^2]+m_2(\tilde{f}_s)\big)ds \Big)\Big(\int_0^t \big(1+C_{\gamma,p}\|f_s\|_{L^{p}}\big)ds\Big)^{\frac{|\gamma|}{2}}\\
\le& \ell^{\gamma/2}\Big(1+ 2m_2(\tilde{f}_0)T \Big)\Big(1+ \int_0^t \big(1+C_{\gamma,p}\|f_s\|_{L^{p}}\big)ds\Big) \le C(T,\tilde{f}_0, f)\ell^{-\kappa(\gamma)}.
\end{align*}
For $I_t^{4,\ell}$, we use the triple H\"{o}lder inequality to write
\begin{align*}
I_t^{4,\ell}
\le& \bb{E}\Big[\int_0^t |W_{s}^K-X_{s}|^2 ds\Big]^{\frac{1}{2}}
\bb{E}\Big[\int_0^t \int_0^1|X_s-X_s^*(\alpha)|^2d\alpha ds\Big]^{\frac{1+\gamma}{2}}
\bb{E} \Big[\int_0^t\int_0^1\bbd{1}_{\{|W_s^K - W_s^*(\alpha)|^\gamma\ge \ell\}}d\alpha ds\Big]^{\frac{|\gamma|}{2}}.
\end{align*}
Thus $I^{4,\ell}_t\leq C(T, f_0, \tilde{f}_0, f)\ell^{-\kappa(\gamma)}$:
use that $\E[|X_s|^2]=\E_\alpha[|X_s^*|^2]=m_2(\tilde f_0)$, that
$m_2(f_s^K)\le 2m_2(f_s)+2\ca{W}_2^2(f_s,f_s^K)$ %\le C(T, f_0, \tilde{f}_0)
 as before
and treat the last term of the product the same as we study  $I_t^{2, \ell}$.

\vskip1mm

{\it Step 3.}
According to Step 1, we now bound $\Delta_s^K(\alpha)$ by \eqref{rrr1} when $|W_s^K - W_s^*(\alpha)|^\gamma\le \ell$
and by  \eqref{rrr2} when $|W_s^K - W_s^*(\alpha)|^\gamma\ge \ell$:
\begin{align*}
\bb{E}[|W_t^K-X_t|^2] &\le \bb{E}[|W_0-X_0|^2] + C\sum_{i=1}^4 I_t^{i,\ell}  + C K^{1-2/\nu} \bb{E}\Big[\int_0^t \int_0^1|W_s^K - W_s^*(\alpha)|^{2+2\gamma/\nu} d\alpha ds\Big]\\
&\quad~ + C \bb{E}\Big[\int_0^t \int_0^1  (|W_s^K-X_s| ^2 + |W_s^*(\alpha)-X_s^*(\alpha)|^2)\min{\big(|W_s^K - W_s^*(\alpha)|^\gamma, \ell \big)} d\alpha ds\Big].
\end{align*}
It then follows from Step 2 that for all $\ell\ge1$, all $K\ge1$,
\begin{align}
\bb{E}[|W_t^K-X_t|^2]
&\le \cW_2^2 (f_0,\tilde f_0) + C(T, f_0, \tilde{f}_0, f) \ell^{-\kappa(\gamma)}   \label{rara}\\
&+ C K^{1-2/\nu} \bb{E}\Big[\int_0^t \int_0^1|W_s^K - W_s^*(\alpha)|^{2+2\gamma/\nu} d\alpha ds\Big]\notag\\
& + C \bb{E}\Big[\int_0^t \int_0^1  |W_s^K-X_s| ^2 |W_s^K - W_s^*(\alpha)|^\gamma d\alpha ds\Big]\notag\\
&+ C\bb{E}\Big[ \int_0^t \int_0^1 |W_s^*(\alpha)-X_s^*(\alpha)|^2 \min{\big(|W_s^K - W_s^*(\alpha)|^\gamma, \ell \big)} d\alpha ds\Big].\notag
\end{align}
Since $\gamma+\nu>0$, it holds that $2+2\gamma/\nu>0$. As a consequence, like in Step 2,
\[\bb{E}\Big[\int_0^t \int_0^1 |W_s^K-W_s^*(\alpha)|^{2+2\gamma/\nu} d\alpha ds\Big]\le C_T[1+\bb{E}[|W_s^K|^2]+m_2(f_0)]\le C(T,f_0, \tilde{f}_0),\]
which gives  \[\lim_{K\to\infty}K^{1-2/\nu} \bb{E}\Big[\int_0^t \int_0^1|W_s^K - W_s^*(\alpha)|^{2+2\gamma/\nu} d\alpha ds\Big]=0.\]
Moreover, we recall that a.s.
$\int_0^1|W_s^K - W_s^*(\alpha)|^\gamma d\alpha \leq 1+ C_{\gamma,p}\|f_s\|_{L^{p}}$ as in Step 2, whence
\[
\bb{E}\Big[\int_0^t \int_0^1  |W_s^K-X_s| ^2 |W_s^K - W_s^*(\alpha)|^\gamma d\alpha ds\Big]
\le \int_0^t \bb{E}[|W_s^K-X_s| ^2](1+C_{\gamma,p}\|f_s\|_{L^p})ds.
\]
Letting $K\to\infty$, by dominated convergence, we find (recall \eqref{dfjt})
\[
\limsup_{K}\bb{E}\Big[\int_0^t \int_0^1  |W_s^K-X_s| ^2 |W_s^K - W_s^*(\alpha)|^\gamma d\alpha ds\Big]
\le \int_0^t J_s(1+C_{\gamma,p}\|f_s\|_{L^p})ds.
\]
Next, it is obvious that for each $\ell\geq 1$ fixed, for all $s\in[0, T]$, all $\alpha\in[0,1]$,
the function $v\mapsto \min(|v-W_s^*(\alpha)|^\gamma, \ell)$ is bounded and continuous. By
\eqref{prop-law-convergence}, we conclude that
$\lim_{K\to \infty} \bb{E}\big[\min\big(|W_s^K - W_s^*(\alpha)|^\gamma, \ell \big) \big]
=\bb{E}\big[\min\big(|W_s- W_s^*(\alpha)|^\gamma, \ell \big)\big]$ and, by dominated convergence, that,
still for $\ell\geq 1$ fixed,
\begin{align*}
&\lim_{K\to\infty}\bb{E}\Big[ \int_0^t \int_0^1 |W_s^*(\alpha)-X_s^*(\alpha)|^2 \min\big(|W_s^K - W_s^*(\alpha)|^\gamma, \ell \big) d\alpha ds\Big]\\
&= \int_0^t \int_0^1 |W_s^*(\alpha)-X_s^*(\alpha)|^2
\bb{E}\big[\min{\big(|W_s- W_s^*(\alpha)|^\gamma, \ell \big)} \big] d\alpha ds.
\end{align*}
But since $W_s \sim f_s$, we have, for each $\alpha$ fixed,
$\bb{E}[\min{(|W_s- W_s^*(\alpha)|^\gamma, \ell )} ]\leq \int_{\bb{R}^3} |W_s^*(\alpha) - v|^\gamma f_s(dv)\leq
1+ C_{\gamma,p}\|f_s\|_{L^{p}}$ by \eqref{norm-inequality}. Furthermore, we have
$\int_0^1|W_s^*(\alpha)-X_s^*(\alpha)|^2 d\alpha=\E_\alpha[|W_s^*-X_s^*|^2]=\cW_2^2(f_s,\tilde f_s)\le J_s$.
All in all, we have checked that
\begin{align*}
&\lim_{K\to\infty}\bb{E}\Big[ \int_0^t \int_0^1 |W_s^*(\alpha)-X_s^*(\alpha)|^2 \min\big(|W_s^K - W_s^*(\alpha)|^\gamma, \ell \big) d\alpha ds\Big] \leq C \int_0^t J_s(1+\|f_s\|_{L^p})ds.
\end{align*}
Gathering all the previous estimates to let $K\to \infty$ in \eqref{rara}: for each $\ell\geq 1$
fixed,
$$
J_t \leq \cW_2^2(f_0,\tilde f_0) + C(T, f_0, \tilde{f}_0, f) \ell^{-\kappa(\gamma)} +C\int_0^t J_s(1+\|f_s\|_{L^p})ds.
$$
Letting now $\ell\to \infty$ and using the  Gr\"{o}nwall lemma, we find
\[J_t \le\ca{W}_2^2(f_0, \tilde{f}_0)
\exp{\left(C_{\gamma,p}\int_0^t\big(1+ \|f_s\|_{L^{p}}\big)ds\right)}.\]
Since $\ca{W}_2^2(f_t, \tilde{f}_t)\le J_t$, this completes the proof.
\epf

It remains to prove Proposition \ref{prop-PDE}. We start with a technical result.

\blem\label{lem-local-condition}
Assume \eqref{con} for some $\gamma\in(-1,0)$, some $\nu\in(0,1)$ with $\gamma+\nu>0$   and  recall that the deviation function $c$ was defined by \eqref{rewrite}. Consider $f\in \ca{P}_2({\bb{R}^3})$
and $\phi_\epsilon(x)=(2\pi\epsilon)^{-3/2} e^{-|x|^2/(2\epsilon)}$. Set  $f^\epsilon(w)=(f * \phi_\epsilon)(w)$.
\benu[label=\emph{(\roman*)}]
\item There exists a constant $C>0$ such that for all $x\in\bb{R}^3$,  all $\epsilon\in(0,1)$,  \begin{align*}
\int_{\bb{R}^3}  \int_{\bb{R}^3} \int_0^\infty \int_0^{2\pi} |c(v,v_*,z, \varphi)| \frac{\phi_\epsilon(v-x)} {f^\epsilon(x)}   d\varphi dz f(dv) f(dv_*)\le C \left(1+\sqrt{m_2(f)}+|x| \right),
\end{align*}
\item For all $\epsilon\in(0,1)$,  all $R>0$, there is a constant $C_{R, \epsilon}>0$ (depending only on
$m_2(f)$) such that for all $x,y\in B(0, R)$,
\begin{align*}
\int_{\bb{R}^3}  \int_{\bb{R}^3} \int_0^\infty \int_0^{2\pi} |c(v,v_*,z, \varphi)| \left|\frac{\phi_\epsilon(v-x)}{f^\epsilon(x)} - \frac{ \phi_\epsilon(v-y)}{f^\epsilon(y)}\right| d\varphi dz f(dv) f(dv_*)\le C_{R, \epsilon}|x-y|.
\end{align*}
\eenu
\elem
\bpf
We start with (i) and  set $I_\e(x)=\int_{\bb{R}^3}  \int_{\bb{R}^3} \int_0^\infty \int_0^{2\pi}
|c(v,v_*,z, \varphi)| \frac{\phi_\epsilon(v-x)} {f^\epsilon(x)}   d\varphi dz f(dv) f(dv_*)$.
Using \eqref{num} and \eqref{con1}, we see that $|c(v,v_*,z, \varphi)|\leq G(z/|v-v_*|^\gamma)|v-v_*|
\leq C(1+z/|v-v_*|^\gamma)^{-1/\nu}|v-v_*|$. Hence
\begin{align*}
I_\e(x)\leq& C \int_{\bb{R}^3}  \int_{\bb{R}^3} \int_0^\infty (1+z/|v-v_*|^\gamma)^{-1/\nu}|v-v_*|\frac{\phi_\epsilon(v-x)}{f^\epsilon(x)} dz f(dv) f(dv_*)\\
 =&C  \int_{\bb{R}^3}  \int_{\bb{R}^3} |v-v_*|^{1+\gamma}\frac{\phi_\epsilon(v-x)}{f^\epsilon(x)} f(dv) f(dv_*).
\end{align*}
Using now that $|v-v_*|^{1+\gamma} \le 1+|v|+|v_*|$, we find
\begin{align*}
I_\e(x)\leq & C \int_{\bb{R}^3}  \int_{\bb{R}^3} (1+|v|+|v_*|)\frac{\phi_\epsilon(v-x)}{f^\epsilon(x)} f(dv) f(dv_*)
\leq
C \Big(1+\sqrt{m_2(f)}+ \frac{\int_{\bb{R}^3}|v|\phi_\epsilon(v-x) f(dv)}{f^\epsilon(x)} \Big).
\end{align*}
To conclude the proof of (i), it remains to study $J_\e(x)= (f^\epsilon(x))^{-1}\int_{\bb{R}^3}|v|\phi_\epsilon(v-x)f(dv)$.
We introduce  $L:=\sqrt{2m_2(f)}$, for which $f(B(0,L))\ge 1/2$ (because $f(B(0,L)^c)\leq m_2(f)/L^2$).
Using that $\{v \in \bb{R}^3:|v|\leq 2|x|+L\}\cup\{v \in \bb{R}^3:|v-x|\geq|x|+L\}=\bb{R}^3$, we write
\[
J_\epsilon(x)= \frac{\int_{\bb{R}^3}|v|\phi_\epsilon(v-x) f(dv)}{\int_{\bb{R}^3}\phi_\epsilon(v-x) f(dv)} \le 2|x|+L+\frac{\int_{|v-x|\ge |x|+L}|v|\phi_\epsilon(v-x) f(dv)}{\int_{|v-x|\le |x|+L} \phi_\epsilon(v-x) f(dv)}.
 \]
Since $\phi_\epsilon$ is radial and decreasing,
\[\int_{|v-x|\ge |x|+L}|v|\phi_\epsilon(v-x) f(dv) \le \phi_\epsilon(|x|+L)\sqrt{m_2(f)}\]
and
\[\int_{|v-x|\le |x|+L} \phi_\epsilon(v-x) f(dv) \ge \phi_\epsilon(|x|+L) f(B(x, |x|+L)) \ge \phi_\epsilon(|x|+L)/2\]
owing to  the fact that $B(0,L)\subset B(x, |x|+L)$. Hence, $J_\e(x)\leq 2|x|+L+2\sqrt{m_2(f)}
\leq 2|x|+4\sqrt{m_2(f)}$
and this completes the proof of (i).

\vskip1mm

For point (ii), we set $\Delta_\e(x,y)=\int_{\bb{R}^3}  \int_{\bb{R}^3} \int_0^\infty \int_0^{2\pi} |c(v,v_*,z, \varphi)|
|F_\e(x,v)-F_\e(y,v)| d\varphi dz f(dv) f(dv_*)$, where
$F_{\epsilon}(v,x):= (f^\epsilon(x))^{-1} \phi_\epsilon(v-x)$.
Exactly as in point (i), we start with
\begin{align*}
\Delta_\e(x,y)
&\le C\int_{\bb{R}^3}  \int_{\bb{R}^3} |v-v_*|^{1+\gamma} |F_{\epsilon}(v,x) - F_{\epsilon}(v,y)| f(dv) f(dv_*)\\
&\le C\int_{\bb{R}^3}(1+\sqrt{m_2(f)} + |v|) |F_{\epsilon}(v,x) - F_{\epsilon}(v,y)| f(dv)\\
&\le C|x-y| \int_{\bb{R}^3}(1+\sqrt{m_2(f)} + |v|)\Big(\sup_{a\in B(0, R)} | \triangledown_x F_{\epsilon}(v,a)| \Big) f(dv)
\end{align*}
for all $x,y \in B(0,R)$.
But we have
\begin{align}\label{gg}
 \triangledown_x F_{\epsilon}(v,a)  =\frac{1}{\epsilon}\frac{\phi_\epsilon(v-a)\int_{\bb{R}^3} (v-u)\phi_\epsilon(u-a) f(du)}{(f^\epsilon(a))^2}.
\end{align}
Indeed, recalling  that $\phi_\epsilon(x)=(2\pi\epsilon)^{-3/2} e^{-|x|^2/(2\epsilon)}$,  we observe that
\[ \triangledown_x \phi_\epsilon(v-x)= \frac{1}{\epsilon} (v-x) \phi_\epsilon(v-x) \  \text {and}\  \triangledown_x f^\epsilon(x)=\frac{1}{\epsilon} \int_{\bb{R}^3}\phi_\epsilon(u-x)(u-x) f(du). \]
Since $F_{\epsilon}(v,a):= (f^\epsilon(a))^{-1}\phi_\epsilon(v-a)$, we have
\begin{align*}
 \triangledown_x F_{\epsilon}(v,a)
 &=\frac{ \triangledown_x \phi_\epsilon(v-a) f^\epsilon(a) - \phi_\epsilon(v-a) \triangledown_x f^\epsilon(a)}{(f^\epsilon(a))^2}\\
 &=\frac{\phi_\epsilon(v-a)}{\epsilon}  \frac{(v-a) f^\epsilon(a) -  \int_{\bb{R}^3}\phi_\epsilon(u-a)(u-a) f(du)}{(f^\epsilon(a))^2}\\
 &=\frac{\phi_\epsilon(v-a)}{\epsilon}  \frac{ \int_{\bb{R}^3}\phi_\epsilon(u-a)(v-a) f(du) -  \int_{\bb{R}^3}\phi_\epsilon(u-a)(u-a) f(du)}{(f^\epsilon(a))^2},
\end{align*}
whence \eqref{gg}. Using now that $J_\e(a)= (f^\epsilon(a))^{-1}\int_{\bb{R}^3} |u| \phi_\epsilon(u-a) f(du)
\le 2|a|+ 4\sqrt{m_2(f)}$ as proved in (i),
\begin{align*}
| \triangledown_x F_{\epsilon}(v,a)|
&\le \frac{1}{\epsilon}\frac{\phi_\epsilon(v-a)}{f^\epsilon(a)}\frac{\int_{\bb{R}^3} (|v|+|u|) \phi_\epsilon(u-a) f(du)}{f^\epsilon(a)}
\le \frac{1}{\epsilon}\frac{\phi_\epsilon(v-a)}{f^\epsilon(a)} \left(|v| +2|a|+ 4\sqrt{m_2(f)}\right).
\end{align*}
But we know that $\phi_\epsilon(x) \le (2\pi\epsilon)^{-3/2}$  and that
\begin{align*}
f^\epsilon(a) \ge \int_{|v-a|\le |a|+L} \phi_\epsilon(v-a) f(dv) \ge \phi_\epsilon(|a|+L) f(B(a, |a|+L)) \ge \phi_\epsilon(|a|+L)/2
\end{align*}
since $B(0,L)\subset B(a, |a|+L)$. Hence,
\[ \sup_{a\in B(0, R)} | \triangledown_x F_{ \epsilon}(v,a)| \le \frac{2}{\epsilon} \ e^{(R+L)^2/(2\epsilon)} \left(|v| +2R+ 4\sqrt{m_2(f)}\right).\]
Consequently, for all $x,y \in B(0,R)$,
\begin{align*}
\Delta_\e(x,y)
&\le \frac{2C}{\epsilon} \  e^{(R+L)^2/(2\epsilon)} |x-y| \int_{\bb{R}^3}\left(1+\sqrt{m_2(f)} + |v|\right) \left(|v| +2R+ 4\sqrt{m_2(f)}\right) f(dv)\leq C_{R,\e} |x-y|,
\end{align*}
where $C_{R,\e}$ depends only on $R,\e$ and  $m_2(f)$ (recall that $L:=\sqrt{2m_2(f)}$).
\epf

Finally, we end the section with the

\bpf[Proof of Proposition \ref{prop-PDE}]
We consider any given weak solution $(\tilde f_t)_{t\geq 0}\in L^\infty([0,\infty),\ca{P}_2(\bb{R}^3))$ to \eqref{Bol}
and we write the proof in several steps.
\vskip1mm

{\it Step 1.} We introduce $\phi_\epsilon(x)=(2\pi\epsilon)^{-3/2} e^{-|x|^2/(2\epsilon)}$
and  $\tilde f_t^\epsilon(w)=(\tilde f_t * \phi_\epsilon)(w)$.
For each $t\ge 0$, we see that  $\tilde f_t^\epsilon$ is a positive smooth function.  We claim that
for any  $\psi \in Lip(\bb{R}^3)$,
$$
\frac{\partial}{\partial t} \int_{\bb{R}^3} \psi(w) \tilde f_t^\epsilon(dw)= \int_{\bb{R}^3} \tilde f_t^\epsilon (dw)
\ca{\tilde{A}}_{t,\epsilon} \psi(w),
$$
where
\beq\label{pro-generator}
\ca{\tilde{A}}_{t,\epsilon} \psi(w)= \int_{\bb{R}^3} \int_{\bb{R}^3} \int_0^\infty\int_0^{2\pi}
[\psi(w+c(v, v_*, z, \varphi)) - \psi(w)] \frac{ \phi_\epsilon(v-w)}{\tilde f_t^{\epsilon}(w)}
d\varphi dz \tilde f_t(dv_*)\tilde f_t(dv).
\eeq
Indeed,  $\tilde f_t^\epsilon(w)=\int_{\bb{R}^3} \phi_\epsilon(v-w) \tilde f_t(dv)$ since $ \phi_\epsilon(x)$ is even.
According to  \eqref{weak} and \eqref{opera}, we have
\begin{align*}
\frac{\partial}{\partial t} \tilde f_t^\epsilon(w)
&=\int_{\bb{R}^3} \int_{\bb{R}^3} \int_0^\infty \int_0^{2\pi} [\phi_\epsilon(v-w+c(v,v_*,z,\varphi))-\phi_\epsilon(v-w)] d\varphi dz \tilde f_t(dv_*) \tilde f_t(dv)\\
&= \int_{\bb{R}^3} \int_0^K \int_0^{2\pi} \left[\int_{\bb{R}^3}\phi_\epsilon(v-w+c(v,v_*,z,\varphi)) \tilde f_t(dv) - \tilde f_t^\epsilon(w) \right] d\varphi dz\tilde  f_t(dv_*)\\
&\quad~~ + \int_{\bb{R}^3} \int_{\bb{R}^3} \int_K^\infty \int_0^{2\pi} [\phi_\epsilon(v-w+c(v,v_*,z,\varphi))-\phi_\epsilon(v-w)] d\varphi dz \tilde f_t(dv_*) \tilde f_t(dv)
\end{align*}
for any $K\geq 1$. We thus have,  for any  $\psi \in Lip(\bb{R}^3)$,
\begin{align*}
\frac{\partial}{\partial t} \int_{\bb{R}^3}\psi(w) &\tilde f_t^\epsilon(dw)
=\int_{\bb{R}^3} \int_{\bb{R}^3} \int_0^K \int_0^{2\pi} \int_{\bb{R}^3}\phi_\epsilon(v-w+c(v,v_*,z,\varphi))\psi(w)
\tilde f_t(dv) d\varphi dz \tilde f_t(dv_*) dw\\
&- \int_{\bb{R}^3} \int_{\bb{R}^3} \int_0^K \int_0^{2\pi} \psi(w)\tilde f^\e_t(w)  d\varphi dz \tilde f_t(dv_*)dw\\
&+ \int_{\bb{R}^3} \int_{\bb{R}^3} \int_{\bb{R}^3} \int_K^\infty \int_0^{2\pi} [\phi_\epsilon(v-w+c(v,v_*,z,\varphi))-\phi_\epsilon(v-w)] \psi(w) d\varphi dz \tilde f_t(dv_*) \tilde f_t(dv) dw.
\end{align*}
Using the change of variables $w-c(v,v_*,z,\varphi) \mapsto w$, we see that the first integral of the RHS equals
$$
\int_{\bb{R}^3} \int_{\bb{R}^3} \int_0^K \int_0^{2\pi} \int_{\bb{R}^3}\phi_\epsilon(v-w) \psi(w
+ c(v,v_*,z,\varphi))\tilde f_t(dv)d\varphi dz \tilde f_t(dv_*) dw.
$$
Consequently,
\begin{align*}
&\frac{\partial}{\partial t} \int_{\bb{R}^3}\psi(w) \tilde f_t^\epsilon(dw)\\
&= \int_{\bb{R}^3}  \int_{\bb{R}^3} \int_0^K \int_0^{2\pi}  \left[\int_{\bb{R}^3} \psi(w+c(v,v_*,z,\varphi))\frac{\phi_\epsilon(v-w) }{\tilde f_t^\epsilon(w)} \tilde f_t(dv) -\psi(w) \right] \tilde f_t^\epsilon(w)  d\varphi dz \tilde f_t(dv_*) dw \\
&\quad + \int_{\bb{R}^3} \int_{\bb{R}^3} \int_{\bb{R}^3} \int_K^\infty \int_0^{2\pi} [\phi_\epsilon(v-w+c(v,v_*,z,\varphi))-\phi_\epsilon(v-w)] \psi(w) d\varphi dz \tilde f_t(dv_*) \tilde f_t(dv) dw\\
&=\int_{\bb{R}^3}  \int_{\bb{R}^3} \int_0^K \int_0^{2\pi}  \int_{\bb{R}^3} \left[\psi(w+c(v,v_*,z,\varphi)) -\psi(w) \right] \frac{\phi_\epsilon(v-w) }{\tilde f_t^\epsilon(w)}  \tilde f_t(dv) d\varphi dz \tilde f_t(dv_*) \tilde f_t^\epsilon(dw)\\
&\quad + \int_{\bb{R}^3} \int_{\bb{R}^3} \int_{\bb{R}^3} \int_K^\infty \int_0^{2\pi} [\phi_\epsilon(v-w+c(v,v_*,z,\varphi))-\phi_\epsilon(v-w)] \psi(w) d\varphi dz \tilde f_t(dv_*) \tilde f_t(dv) dw.
\end{align*}
Letting $K$ increase to infinity, one easily ends the step.

\vskip1mm {\it Step 2.} We set $F_{t,\epsilon}(v,x)= (\tilde f_t^\epsilon(x))^{-1} \phi_\epsilon(v-x)$.
For a given  $X_0^\epsilon$ with law $\tilde f_0^\epsilon$, and  a given independent Poisson measure $N(ds, dv, dv_*, dz, d\varphi, du)$ on $[0,\infty) \times \bb{R}^3 \times \bb{R}^3 \times [0, \infty) \times [0, 2\pi) \times [0, \infty)$ with intensity $ds \tilde f_s(dv) \tilde f_s(dv_*) dz d\varphi du$, there  exists a  pathwise unique solution to
\beq\label{pro-convolution-process}
X_t^\epsilon=X_0^\epsilon + \int_0^t \int_{\bb{R}^3} \int_{\bb{R}^3} \int_0^\infty \int_0^{2\pi} \int_0^\infty c(v, v_*, z, \varphi)\bbd{1}_{\left \{u\le F_{s,\epsilon}(v, X_{s-}^\epsilon)\right \}} N(ds, dv, dv_*, dz, d\varphi, du).
\eeq
This classically follows from Lemma \ref{lem-local-condition}, which precisely tells us that
the coefficients of this equation satisfy some
{\it at most linear growth condition} (point (i)) and some {\it local Lipschitz condition} (point (ii)).

\vskip1mm

{\it Step 3.} We now prove that $\ca{L}(X_t^\epsilon)=\tilde f_t^\epsilon$ for each $t\ge0$.
We thus introduce $g^\e_t=\ca{L}(X_t^\epsilon)$. By the It\^o formula, we see that for all
$\psi \in Lip({\bb{R}^3})$,
\begin{align*}
&\frac{\partial}{\partial t} \int_{\bb{R}^3} \psi(w) g_t^\epsilon(dw)\\
=& \int_{\bb{R}^3} g_t^\epsilon (dw)
\int_{\bb{R}^3}\int_{\bb{R}^3}\int_0^\infty\int_0^{2\pi} \Big(\psi(w+c(v,v_*,z,\varphi))-\psi(w) \Big)
F_{t,\e}(v,w) d\varphi dz \tilde f_t(dv_*)\tilde f_t(dv)\\
=&  \int_{\bb{R}^3} g_t^\epsilon (dw) \ca{\tilde{A}}_{t,\epsilon} \psi(w).
\end{align*}
Thus $(\tilde f^\e_t)_{t\geq 0}$ and $(g^\e_t)_{t\geq 0}$ satisfy the same equation and we of course have
$g^\e_0=\tilde f^\e_0$ by construction. The following uniqueness result allows us to conclude the step:
for any $\mu_0 \in \ca{P}_2({\bb{R}^3})$, there exists at most one family
$(\mu_t)\in L_{loc}^\infty\big([0,\infty), \ca{P}_2({\bb{R}^3})\big)$ such that for any
$\psi \in Lip(\bb{R}^3)$,  any $t\ge0$,
\beq\label{linear-SIE}
\int_{\bb{R}^3} \psi(w) \mu_t(dw)= \int_{\bb{R}^3} \psi(w) \mu_0 (dw)
+ \int_0^t ds\int_{\bb{R}^3} \mu_s(dw) \tilde{\ca{A}}_{s,\epsilon}\psi(w).
\eeq
This must be classical (as well as Step 2 is),  but we find no precise reference
and thus make use of martingale problems.
A c\`{a}dl\`{a}g adapted $\bb{R}^3$-valued process $(Z_t)_{t\ge0}$ on some filtered probability space
$(\Omega, \ca{F}, \ca{F}_t, \bb{P})$ is said to solve the martingale problem $MP(\tilde{\ca{A}}_{t,\epsilon}, \mu_0,
Lip(\bb{R}^3))$ if $\bb{P}\circ Z_0=\mu_0$ and if for all $\psi \in Lip(\bb{R}^3)$,
$(M_t^{\psi,\epsilon})_{t\ge0}$ is a $(\Omega, \ca{F}, \ca{F}_t, \bb{P})$-martingale, where
\[M_t^{\psi, \e}=\psi(Z_t) - \int_0^t  \tilde{\ca{A}}_{s,\epsilon}\psi(Z_s)ds.\]
According to \cite[Theorem 5.2]{MR1245309} (see also \cite[Remark 3.1, Theorem 5.1]{MR1245309}
and \cite[Theorem B.1]{MR1042343}), it suffices to check the following points to conclude the uniqueness
for \eqref{linear-SIE}.

\benu[label=(\roman*)]

\item there exists a countable family $(\psi_k)_{k\ge1}\subset Lip(\bb{R}^3)$ such that for all $t\ge0$,
the closure (for the bounded pointwise convergence) of $\{(\psi_k,  \tilde{\ca{A}}_{t,\epsilon}\psi_k), k\ge 1\}$
contains $\{(\psi,  \tilde{\ca{A}}_{t,\epsilon}\psi), \psi\in Lip(\bb{R}^3)\}$,

\item for each $w_0\in\bb{R}^3$, there exists a solution to $MP( \tilde{\ca{A}}_{t,\epsilon}, \delta_{w_0}, Lip(\bb{R}^3))$,

\item for each $w_0\in\bb{R}^3$, uniqueness (in law) holds for $MP( \tilde{\ca{A}}_{t,\epsilon}, \delta_{w_0}, Lip(\bb{R}^3))$.
\eenu

The fact that \eqref{pro-convolution-process} has a pathwise unique solution proved in Step 2
(there we can of course replace $X_0^\e$ by any deterministic point $w_0\in\bb{R}^3$) immediately
implies (ii) and (iii). Point (i) is very easy (recall that $\e>0$ is fixed here).

\vskip1mm {\it Step 4.} In this step, we check that the family $((X_t^\epsilon)_{t\ge0})_{\epsilon>0}$ is tight in $\bb{D}([0,\infty), \bb{R}^3)$. To do this, we use the Aldous criterion \cite{MR0474446}, see also \cite[p 321]{MR1943877}, i.e. it suffices to prove that for all $T>0$,
\beq\label{aldous}
\sup_{\epsilon \in(0,1)}\bb{E}\big [\sup_{[0,T]}|X_t^\epsilon| \big ]< \infty, \quad
\lim_{\delta\to 0}\sup_{\epsilon \in(0,1)}\sup_{S, S^\prime\in \ca{S}_T(\delta)}\bb{E}\big [|X_{S^\prime}^\epsilon-X_S^\epsilon| \big]=0,
\eeq
where $\ca{S}_T(\delta)$ is the set containing all pairs of stopping times $(S, S^\prime)$ satisfying $0\le S\le S^\prime \le S+\delta\le T$. 

\vskip1mm

First,  since $X^\e_t\sim \tilde f^\e_t= \tilde f_t \star \phi_\e$,
we have $\E[|X^\e_t|^2]\leq 2 (m_2(\tilde f_t) + 3\e) \leq 2m_2(\tilde f_0) + 6$.
Thus for any $T>0$, using Lemma \ref{lem-local-condition}-(i),
\begin{align*}
\bb{E} \Big[\sup_{[0,T]}|X_t^\epsilon| \Big]
&\le \bb{E} \big[|X_0^\epsilon| \big]
+ \bb{E} \Big[\int_0^T\int_{\bb{R}^3}\int_{\bb{R}^3}\int_0^\infty\int_0^{2\pi}|c(v,v_*,z, \varphi)|
 \frac{\phi_\epsilon(v-X_s^{\epsilon})}{\tilde f_s^\epsilon(X_s^{\epsilon})}  d\varphi dz \tilde f_s(dv) \tilde f_s(dv_*) ds \Big]\\
&\le  \bb{E} \big[|X_0^\epsilon| \big] + C \bb{E}\left[\int_0^T (1+|X_s^\epsilon|) \, ds \right] \leq C_T.
\end{align*}
Furthermore,  for any $T>0$, $\delta>0$ and $(S, S^\prime)\in\ca{S}_T(\delta)$, using again
Lemma \ref{lem-local-condition}-(i),
\begin{align*}
\bb{E} \big [|X_{S^\prime}^\epsilon-X_S^\epsilon| \big]
&\le \bb{E}\left[\int_S^{S+\delta}\int_{\bb{R}^3}\int_{\bb{R}^3}\int_0^\infty\int_0^{2\pi}|c(v,v_*,z, \varphi)| \frac{\phi_\epsilon(v-X_s^{\epsilon})}{\tilde f_s^\epsilon(X_s^{\epsilon})}  d\varphi dz \tilde f_s(dv) \tilde f_s(dv_*) ds \right]\\
&\le C\bb{E}\left[\int_S^{S+\delta}(1+|X_s^\epsilon|) ds \right] \le C\bb{E}\left[\delta \sup_{[0,T]}(1+|X_s^\epsilon |)  \right] \le C_T \delta.
\end{align*}
Hence \eqref{aldous} holds true and this completes the step.

\vskip1mm 

{\it Step 5.} We thus can find some $(X_t)_{t\ge0}$ which is the limit in law (for the Skorokhod topology)
of a sequence $(X_t^{\epsilon_n})_{t\ge0}$ with $\epsilon_n\searrow 0$. Since $\ca{L}(X_t^{\epsilon_n})=
\tilde f_t^{\epsilon_n}$ by Step 3 and since $\tilde f_t^{\epsilon_n} \to \tilde f_t$ by definition,
we have $\ca{L}(X_t)=\tilde f_t$ for each $t\ge0$. It only remains to show that $(X_t)_{t\ge0}$ is a (weak) solution
to \eqref{SDE}. Using the theory of martingale problems, it classically suffices to prove that
for any $\psi \in C^1_b(\bb{R}^3)$, the process $\psi(X_t)-\psi(X_0) - \int_0^t \ca{B}_s\psi(X_s)ds$ is a martingale,
where
$$
\ca{B}_t\psi(x)=\int_0^1\int_0^\infty\int_0^{2\pi}\Big(\psi(x+c(x,X_t^*(\alpha),z,\varphi))-\psi(x)\Big)
d\varphi dz d\alpha.
$$
But since $\cL_\alpha(X_t^*)=\tilde f_t$, this rewrites (recall \eqref{opera})
$$
\ca{B}_t\psi(x)=\int_{\bb{R}^3}\int_0^\infty\int_0^{2\pi}\Big(\psi(x+c(x,v_*,z,\varphi)-\psi(x)\Big)
d\varphi dz \tilde f_t(dv_*)= \int_{\bb{R}^3} \ca{A} \psi (x,v_*) \tilde f_t(dv_*).
$$
We thus have to prove that for any $0\le s_1\le ... \le s_k \le s\le t\le T$,
any $\psi_1,..., \psi_k\in C_b^1(\bb{R}^3)$, and any $\psi\in C_b^1(\bb{R}^3)$,
$$
\E[\ca{F}(X)]=0,
$$
where $\ca{F}: \bb{D}([0,\infty), \bb{R}^3)\mapsto \bb{R}$ is defined by
\[\ca{F}(\lambda)=\Big(\prod_{i=1}^{k}\psi_i(\lambda_{s_i})\Big) \Big(\psi(\lambda_t)-\psi(\lambda_s)-
\int_s^t \ca{B}_r\psi(\lambda_r) dr\Big).\]
We of course start from $\bb{E}[\ca{F}_{\epsilon_n}(X^{\epsilon_n})]=0$, where, recalling
\eqref{pro-generator},
\[\ca{F}_\epsilon(\lambda)=\Big(\prod_{i=1}^{k}\psi_i(\lambda_{s_i})\Big) \Big(\psi(\lambda_t)-\psi(\lambda_s)-\int_s^t \tilde{\ca{A}}_{r,\e} \psi(\lambda_r) dr\Big).\]
We then write
\begin{align*}
\Big| \bb{E}[\ca{F}(X)] \Big| \le \Big |\bb{E}[\ca{F}(X)]-\bb{E}[\ca{F}(X^{\epsilon_n})] \Big |
+\Big |\bb{E}[\ca{F}(X^{\epsilon_n})]- \bb{E}[\ca{F}_{\epsilon_n}(X^{\epsilon_n})]\Big |.
\end{align*}

On the one hand,
we know from \cite[Lemma 3.3]{MR3313757} that $(x,v_*)\mapsto \ca{A} \psi(x,v_*)$ is continuous on
$\bb{R}^3\times\bb{R}^3$ and bounded by $C\, |x-v_*|^{\gamma+1}$.
We thus  easily deduce that $\ca{F}$ is continuous at each $\lambda \in\bb{D}([0,\infty), \bb{R}^3)$
which does not jump at $s_1,..., s_k, s, t$ (this is a.s. the case of $X\in\bb{D}([0,\infty), \bb{R}^3)$
because it has no deterministic time jump by the Aldous criterion). We also
deduce that $|\ca{F} (\lambda)|\leq C(1+\int_0^t \int_{\bb{R}^3} |\lambda_r - v_*|^{\gamma+1} \tilde f_r(dv_*)dr)$.
Using that $0<\gamma+1<1$, that $\sup_{\e\in(0,1)}\E[\sup_{[0,T]} |X^{\e}_t|]<\infty$ by Step 4
and recalling that $X^{\e_n}$ goes in law to $X$, we easily conclude that
$|\bb{E}[\ca{F}(X)]-\bb{E}[\ca{F}(X^{\epsilon_n})] |$ tends to $0$ as $n\to \infty$.

\vskip1mm

On the other hand, since $|\ca{F}(\lambda)-\ca{F}_\e(\lambda)|\leq C
|\int_s^t (\ca{B}_r\psi(\lambda_r)-\tilde{\ca{A}}_{r,\e}\psi(\lambda_r))dr|$
and  $X^\e_r\sim \tilde f^\e_r$,
\begin{align*}
&\Big |\bb{E}[\ca{F}(X^{\epsilon_n})]- \bb{E}[\ca{F}_{\epsilon_n}(X^{\epsilon_n})]\Big |\\
\le& C  \int_s^t \E\Big[ \Big |
\int_{\bb{R}^3}\int_0^\infty \int_0^{2\pi} \int_{\bb{R} ^3}\psi(X_r^{\epsilon_n}+c(v, v_*, z, \varphi))
\Big[\frac{ \phi_{\epsilon_n}(v-X_r^{\epsilon_n})}{\tilde f_r^{\epsilon_n}(X_r^{\epsilon_n})}\tilde f_r(dv) -
\dirac_{X_r^{\epsilon_n}}(dv)\Big]
d\varphi dz \tilde f_r(dv_*) \Big|\Big] dr \\
=&C \int_s^t \Big| \int_{\bb{R}^3}\int_0^\infty \int_0^{2\pi} \int_{\bb{R}^3} \int_{\bb{R} ^3}\psi(w+c(v, v_*, z, \varphi))
\Big[ \phi_{\epsilon_n}(v-w)\tilde f_r(dv)- \tilde f_r^{\epsilon_n}(w) \dirac_{w}(dv)\Big]
dw d\varphi dz \tilde f_r(dv_*) \Big| dr.
\end{align*}
But we can write $\int_{\bb{R}^3}\int_{\bb{R}^3}\psi(w+c(v,v_*,z,\varphi))\tilde f_r^{\epsilon_n}(w) \dirac_{w}(dv)dw
=\int_{\bb{R}^3}\psi(w+c(w,v_*,z,\varphi))\tilde f_r^{\epsilon_n}(w)dw=\int_{\bb{R}^3}\int_{\bb{R}^3}
\psi(w+c(w,v_*,z,\varphi))\phi_{\e_n}(v-w)\tilde f_r(dv)dw$, so that
\begin{align*}
&\Big |\bb{E}[\ca{F}(X^{\epsilon_n})]- \bb{E}[\ca{F}_{\epsilon_n}(X^{\epsilon_n})]\Big |\\
\leq & C\int_s^t \Big|\int_{\bb{R}^3}\int_0^\infty \int_0^{2\pi} \int_{\bb{R}^3}\int_{\bb{R} ^3} \Big[ \psi(w+c(v, v_*, z, \varphi)) - \psi(w+c(w, v_*, z, \varphi))\Big]\\
&\hskip7cm  \phi_{\epsilon_n}(v-w)\tilde f_r(dv) dw d\varphi dz \tilde f_r(dv_*) \Big| dr\\
=& C\int_s^t \Big|\int_{\bb{R}^3}\int_0^\infty \int_0^{2\pi} \int_{\bb{R}^3}\int_{\bb{R} ^3} \Big[ \psi(w+c(v, v_*, z, \varphi+\varphi_0(v-v_*, w-v_*))) \\
&\hskip 4cm  - \psi(w+c(w, v_*, z, \varphi))\Big]  \phi_{\epsilon_n}(v-w)\tilde f_r(dv) dw d\varphi dz \tilde f_r(dv_*)\Big | dr.
\end{align*}
The last equality uses the $2\pi$-periodicity of $c$.
We now put
$$
R_n(v,v_*,z,\varphi):=\int_{\bb{R} ^3} \Big[ \psi(w+c(v, v_*, z, \varphi+\varphi_0(v-v_*, w-v_*))) - \psi(w+c(w, v_*, z, \varphi))\Big] \phi_{\epsilon_n}(v-w) dw,
$$
and show the following two things:
\benu[label=(\alph*)]
\item for all $v, v_*\in\bb{R}^3$, all $z\in[0,\infty)$ and $\varphi\in[0,2\pi)$,
$\lim_{n\to\infty}R_n(v,v_*,z,\varphi)=0$;
\item there is a constant $C>0$ such that for all $n\geq 1$, all
$v, v_*\in\bb{R}^3$, all $z\in[0,\infty)$ and $\varphi\in[0,2\pi)$,
$$
|R_n(v,v_*,z,\varphi)|\le  C\big(1+|v-v_*|\big)(1+z)^{-1/\nu},
$$
which belongs to
$L^1([0,T]\times\bb{R}^3\times\bb{R}^3\times[0,\infty)\times[0,2\pi),dr \tilde f_r(dv_*)\tilde f_r(dv)
dz d\varphi)$ because $(\tilde f_t)_{t\geq 0} \in L^\infty([0,T],\ca{P}_2(\bb{R}^3))$ by assumption.
\eenu

By dominated convergence, we will deduce that
$\lim_{n\rightarrow \infty}\Big |\bb{E}[\ca{F}(X^{\epsilon_n})]- \bb{E}[\ca{F}_{\epsilon_n}(X^{\epsilon_n})]\Big |=0$
and this will conclude the proof.

\vskip1mm

We first study (a). Since $\psi\in C^1_b(\bb{R}^3)$, we immediately observe that
\begin{align}\label{startagain}
& \Big| \psi(w+c(v, v_*, z, \varphi+\varphi_0(v-v_*, w-v_*))) - \psi(w+c(w, v_*, z, \varphi))\Big | \\
&\le C_\psi \Big |c(v, v_*, z, \varphi+\varphi_0(v-v_*, w-v_*)) - c(w, v_*, z, \varphi)\Big|.\notag
\end{align}
Recalling  that
\[c(v,v_*,z, \varphi)=-\frac{1-\cos G(z/|v-v_*|^\gamma)}{2}(v-v_*)+\frac{ \sin G(z/|v-v_*|^\gamma)) }{2} \Gamma(v-v_*, \varphi) ,\]
we have
\begin{align*}
& \Big |c(v, v_*, z, \varphi+\varphi_0(v-v_*, w-v_*)) - c(w, v_*, z, \varphi)\Big|\\
\le& \frac{|\cos G(z/|v-v_*|^\gamma) - \cos G(z/|w-v_*|^\gamma)|} {2} |v-v_*|+ \frac{ |1 - \cos G(z/|w-v_*|^\gamma)|}{2} |v-w|\\
& + \frac{|\sin G(z/|v-v_*|^\gamma) - \sin G(z/|w-v_*|^\gamma)|}{2} |\Gamma(v-v_*, \varphi+\varphi_0)| \\
& + \frac{|\sin G(z/|w-v_*|^\gamma)|}{2}|\Gamma(v-v_*, \varphi+\varphi_0) - \Gamma(w-v_*, \varphi)|.
\end{align*}
Using  that $|\Gamma(v-v_*, \varphi+\varphi_0)|=|v-v_*|$ and Lemma \ref{Tanaka},  we obtain
\begin{align*}
& \Big |c(v, v_*, z, \varphi+\varphi_0(v-v_*, w-v_*)) - c(w, v_*, z, \varphi)\Big|\\
\le& C |G(z/|v-v_*|^\gamma) - G(z/|w-v_*|^\gamma)||v-v_*|+ C |v-w|.
\end{align*}
We deduce from \eqref{nota}  that $|G^\prime(z)| =1/\beta(G(z)) \leq C$ by \eqref{con}, whence
\begin{align*}
& \Big |c(v, v_*, z, \varphi+\varphi_0(v-v_*, w-v_*)) - c(w, v_*, z, \varphi)\Big|\\
&\le C z \big | |v-v_*|^{|\gamma|} - |w-v_*|^{|\gamma|}\big| |v-v_*|+ C |v-w|.
\end{align*}
Using  again the inequality  $|x^{\alpha} - y^{\alpha}|\le|x-y|(x\vee y)^{\alpha-1}$ for $\alpha\in(0,1)$, and $x,y\ge0$, we have
\[\big | |v-v_*|^{|\gamma|} - |w-v_*|^{|\gamma|}\big| \le |v-w| |v-v_*|^{|\gamma|-1}.\]
We thus get
\[\Big |c(v, v_*, z, \varphi+\varphi_0(v-v_*, w-v_*)) - c(w, v_*, z, \varphi)\Big|
\le C (z |v-v_*|^{|\gamma|}+1) |v-w|.
\]
Consequently,
\begin{align*}
R_n(v,v_*,z,\varphi)
 \le  &C_\psi  (z |v-v_*|^{|\gamma|}+1)  \int_{\bb{R} ^3} |v-w| \phi_{\epsilon_n}(v-w) dw,
\end{align*}
which clearly tends to $0$ as $n\to \infty$.
This ends the proof of (a).

\vskip1mm

For (b), start again from \eqref{startagain} to write
\begin{align*}
& \Big| \psi(w+c(v, v_*, z, \varphi+\varphi_0(v-v_*, w-v_*))) - \psi(w+c(w, v_*, z, \varphi))\Big | \\
&\le\Big| \psi(w+c(v, v_*, z, \varphi)) -\psi(w) \Big |+ \Big |\psi(w)- \psi(w+c(w, v_*, z, \varphi))\Big | \\
&\le C_\psi ( |c(v, v_*, z, \varphi)| + |c(w, v_*, z, \varphi)|).
\end{align*}
Moreover, since $|c(v, v_*, z, \varphi)|\le G(z/|v-v_*|^\gamma)|v-v_*| \leq C |v-v_*| (1+ |v-v_*|^{|\gamma|}z)^{-1/\nu}$
by \eqref{num} and \eqref{con1}, we observe that
\begin{align*}
R_n(v,v_*,z,\varphi)\leq C |v-v_*| (1+ |v-v_*|^{|\gamma|}z)^{-1/\nu} + C \int_{\bb{R} ^3} |w-v_*|
(1+ |w-v_*|^{|\gamma|}z)^{-1/\nu}\phi_{\epsilon_n}(v-w) dw.
\end{align*}
Since $1+ |v-v_*|^{|\gamma|}z\ge \big(1\wedge |v-v_*|^{|\gamma|}\big)(1+z)$ for $z\in[0,\infty)$,
\begin{align*}
|v-v_*| (1+ |v-v_*|^{|\gamma|}z)^{-1/\nu}  &\le  |v-v_*|(1+z)^{-1/\nu}\big(1\wedge |v-v_*|^{|\gamma|} \big)^{-1/\nu}.
\end{align*}
Using that  $|\gamma|/\nu<1$, we deduce that
$$
|v-v_*| (1+ |v-v_*|^{|\gamma|}z)^{-1/\nu} \le  \big(1+|v-v_*| \big) (1+z)^{-1/\nu}.
$$
As a conclusion,
$$
R_n(v,v_*,z,\varphi)\leq C \Big(1+|v-v_*|+\int_{\bb{R} ^3} |w-v_*| \phi_{\epsilon_n}(v-w) dw \Big) (1+z)^{-1/\nu},
$$
which is easily bounded (recall that $\e_n\in (0,1)$)
by $C(1+|v|+|v_*|)(1+z)^{-1/\nu}$ as desired.
\epf

\section{The coupling}

To get the convergence of the particle system, we construct a suitable coupling between the particle
system with generator $\ca{L}_{N,K}$ defined by \eqref{rewpartopera}
and the realization of the weak solution to \eqref{Bol}, following the ideas of \cite{MR3456347}.

\blem\label{coupling-systems}
Assume \eqref{con} for some $\gamma\in(-1,0)$, $\nu\in(0,1)$ with $\gamma+\nu\in(0,1)$.
Let $N\ge1$ be fixed. Let  $q\ge2$ such that $q>\gamma^2/(\gamma+\nu)$.  Let $f_0\in\ca{P}_q(\bb{R}^3)$
with a finite entropy and let $(f_t)_{t\ge0} \in L^\infty\big([0,\infty), \ca{P}_2({\bb{R}^3})\big)
\cap L^1_{loc}\big([0,\infty), L^{p}({\bb{R}^3})\big)$ (with $p\in(3/(3+\gamma),p_0(\gamma,\nu,q))$)
be the unique weak solution to \eqref{Bol}  given by
Theorem \emph{ \ref{well-posedness}}. Then there exists, on some probability space, a family of  i.i.d.  random
variables $(V_0^i)_{i=1,...,N}$  with common law $f_0$, independent of  a family of i.i.d. Poisson measures
$(M_i(ds, d\alpha, dz, d\varphi))_{i=1,...,N}$ on $[0,\infty)\times[0,1]\times[0,\infty)\times[0,2\pi)$,  with
intensity $ds d\alpha dz d\varphi$, a measurable family $(W^*_t)_{t\geq 0}$ of $\alpha$-random variables with
$\alpha$-law $(f_t)_{t\ge0}$ and $N$ i.i.d. c\`adl\`ag
adapted processes $(W_t^i)_{t\ge0}$  solving, for each $i=1, \cdots, N$,
\beq\label{B-particle}
W_t^i=V_0^i+\int_0^t\int_0^1\int_0^\infty\int_0^{2\pi}c(W_{s-}^{i},W_s^*(\alpha),z,\varphi)M_i(ds,d\alpha,dz,d\varphi).
\eeq
Moreover, $W_t^i \sim f_t$ for all $t\ge0$, all $i=1,\dots,N$.
Also, for all $T>0$,
\beq\label{lem-momen}
\bb{E}\big[\sup_{[0,T]} |W_t^{1}|^q\big]\le C_{T,q}.
\eeq

\elem
\bpf
Except for the moment estimate \eqref{lem-momen}, it suffices to apply Proposition \ref{prop-PDE}.
A simpler proof could be handled here because we deal with the {\it strong} solution
$f \in L^\infty\big([0,\infty), \ca{P}_2({\bb{R}^3})\big)
\cap L^1_{loc}\big([0,\infty), L^{p}({\bb{R}^3})\big)$.
We now prove \eqref{lem-momen}, which is more or less classical. We thus fix $q\geq 2$.
It is clear that
\begin{align*}
\big||v+c(v,v_*,z,\varphi)|^q-|v|^q\big|  \le  C_q  \big( |v|^{q-1}+|c(v,v_*,z,\varphi)|^{q-1} \big) |c(v,v_*,z,\varphi)| .
\end{align*}
Due to \eqref{num} and \eqref{con1},
$|c(v,v_*,z,\varphi)|\le |v-v_*|$,  $|c(v,v_*,z,\varphi)|\le (1+z/|v-v_*|^\gamma)^{-1/\nu}|v-v_*|$, whence
\begin{align}\label{q-increment}
&\int_0^\infty \int_0^{2\pi} \big||v+c(v,v_*,z,\varphi)|^q-|v|^q\big| d\varphi dz \notag \\
&\le C_q  \int_0^\infty \int_0^{2\pi} \big( 1+|v|^{q-1}+|v_*|^{q-1}\big) (1+z/|v-v_*|^\gamma)^{-1/\nu}
|v-v_*| d\varphi dz \notag \\
&= C_q \big( 1+|v|^{q-1}+|v_*|^{q-1}\big) |v-v_*|^{1+\gamma} \notag\\
&\le C_q \big( 1+|v|^{q}+|v_*|^{q}\big),
\end{align}
because $0<1+\gamma < 1$. It then easily
follows  from the It\^{o} formula  and $\ca{L}_\alpha(W_t^*)=f_t=\cL(W_t^1)$ that
\begin{align*}
\bb{E}\big[\sup_{[0,t]}|W_s^{1}|^q\big]
&\le \bb{E}[|V_0^{1}|^q] + C_q \int_0^t\int_0^1\bb{E}\Big[ 1+|W_s^{1}|^q+|W_s^*(\alpha)|^q \Big]  d\alpha ds\\
&\le \bb{E}[|V_0^{1}|^q] + C_q  \int_0^t \Big(1+\bb{E}[\sup_{[0,s]}|W_u^{1}|^q]\Big)  ds.
\end{align*}
We thus conclude \eqref{lem-momen} by the Gr\"{o}nwall lemma.
\epf

Next, let us recall \cite[Lemma 4.3]{MR3456347} below in order to construct our coupling.

\blem\label{coupling}
Consider $(f_t)_{t\geq 0}$ and $(W^*_t)_{t\geq 0}$ introduced in Lemma \emph{\ref{coupling-systems}}
and fix $N\ge1$.
For  ${\bf v}=(v_1, v_2, ...,v_N)\in(\bb{R}^3)^N$,
we introduce the empirical measure $\mu_{{\bf v}}^N:=N^{-1}\sum_{i=1}^{N}\delta_{v_i}$.
Then for all $t\ge0$, all ${\bf v}\in (\bb{R}^3)^N$ and all
${\bf w}\in (\bb{R}^3)^N_{\bullet}$, with $(\bb{R}^3)^N_{\bullet}:=\{{\bf w}\in(\bb{R}^3)^N: w_i\ne w_j \; \forall \,
i\ne j \}$,
there are $\alpha$-random variables
$Z_t^*({\bf w},\alpha)$ and $V_t^*({\bf v},{\bf w},\alpha)$ such that the $\alpha$-law of $(Z_t^*({\bf w},\cdot), V_t^*({\bf v},{\bf w},\cdot))$ is $N^{-1}\sum_{i=1}^{N}\delta_{(w_i, v_i)}$ and $\int_0^1|W_t^*(\alpha)-Z_t^*({\bf w},\alpha)|^2 d\alpha=\ca{W}_2^2(f_t, \mu_{{\bf w}}^N)$.
 \elem

Observe that $\ca{L}_\alpha(Z_t^*({\bf w},\cdot))=\mu_{{\bf w}}^N$ and
$\ca{L}_\alpha(V_t^*({\bf v},{\bf w},\cdot))=\mu_{{\bf v}}^N$.
Owing to technical reasons, we need to introduce some more notations.

\begin{nota}\label{nnn}
We consider an $\alpha$-random variable $Y$ with uniform distribution on $B(0,1)$ (independent of everything else)
and, for $\epsilon\in (0,1)$, $t\ge0$,
$\alpha\in[0,1]$, ${\bf v}\in (\bb{R}^3)^N$ and ${\bf w}\in (\bb{R}^3)^N_{\bullet}$, we set
$W_t^{*,\epsilon}(\alpha)=W_t^*(\alpha)+\epsilon Y(\alpha)$ and $V_t^{*,\epsilon}({\bf v},{\bf w},\alpha)=
V_t^*({\bf v},{\bf w},\alpha)+\epsilon Y(\alpha)$.
It holds that $\ca{L}_\alpha(W_t^{*,\epsilon})=f_t*\psi_\epsilon$ and $\ca{L}_\alpha(V_t^{*,\epsilon}({\bf v},{\bf w},\cdot))
=\mu_{{\bf v}}^N*\psi_\epsilon$, where $\psi_\epsilon(x)=(3/(4\pi\epsilon^3))\bbd{1}_{\{|x|\le \epsilon\}}$.
\end{nota}

At last, we built a suitable realisation for the particle system.

\blem\label{aaaa}
Consider all the objects introduced in Lemmas \emph{\ref{coupling-systems}}-\emph{\ref{coupling}} and Notation \emph{\ref{nnn}}.
Set ${\bf {W}}_{s}=(W_{s}^{1},...,W_{s}^{N})$, which a.s. belongs to $(\bb{R}^3)^N_\bullet$
(because $f_s$ has a density for all $s\geq 0$). Fix $K\geq 1$ and $\e\in (0,1)$.
There is a unique strong solution
$({\bf {V}}_{t})_{t\geq 0}=(V_{t}^{1}, ... , V_{t}^{N})_{t\geq 0}$ to
\beq\label{particle-system}
V_t^{i}=V_0^i+\int_0^t\int_0^1\int_0^\infty\int_0^{2\pi}c_K(V_{s-}^{i}, V_s^*({\bf {V}}_{s-},  {\bf {W}}_{s-}, \alpha), z, \varphi+\varphi_{i,\alpha,s})M_i(ds,d\alpha,dz,d\varphi), ~  i=1,..., N,
\eeq
where $\varphi_{i,\alpha,s}:=\varphi_{i,\alpha,s}^1+\varphi_{i,\alpha,s}^2+\varphi_{i,\alpha,s}^3$ with
\begin{align*}
\varphi_{i,\alpha,s}^1=&\varphi_{0}\big(W_{s-}^{i} - W_s^*(\alpha), W_{s-}^{i} - W_s^{*,\epsilon}(\alpha)\big),\\
\varphi_{i,\alpha,s}^2=&\varphi_{0}\big(W_{s-}^{i} - W_s^{*,\epsilon}(\alpha),
V_{s-}^{i} - V_s^{*,\epsilon}({\bf {V}}_{s-}, {\bf {W}}_{s-}, \alpha)\big),\\
\varphi_{i,\alpha,s}^3=&\varphi_{0}\big( V_{s-}^{i} - V_s^{*,\epsilon}({\bf {V}}_{s-}, {\bf {W}}_{s-}, \alpha),
V_{s-}^{i} - V_s^*({\bf {V}}_{s-}, {\bf {W}}_{s-}, \alpha)\big).
\end{align*}
Moreover, $({\bf {V}}_{t})_{t\geq 0}$ is a Markov process with generator
$\ca{L}_{N,K}$. %We have $\bb{E}\big[|V_t^{1}|^2\big]=m_2(f_0)$ and,
If $f_0\in\ca{P}_q(\bb{R}^3)$
for some $q\ge2$, then $\bb{E}\big[\sup_{[0,T]}|V_t^{1}|^q\big]\le C_{T,q}$
(this last constant not depending on $N,K$ nor $\e\in(0,1)$).
\elem

\bpf
Since $c_K=c\bbd{1}_{\{z\le K\}}$, the Poisson measures involved in \eqref{particle-system} are finite.
Hence the existence and uniqueness results hold for \eqref{particle-system}.
Using that $\ca{L}_\alpha(V_t^*({\bf v},{\bf w},\cdot))=\mu_{{\bf v}}^N$
and the $2\pi$-periodicity of $c_K$ in $\varphi$, one easily checks that $({\bf {V}}_{t})_{t\geq 0}$ is a
Markov process with generator $\ca{L}_{N,K}$: for all bounded measurable function
$\phi:(\bb{R}^3)^N\mapsto \bb{R}$, all $t\geq 0$, a.s.,
\begin{align*}
&\sum_{i=1}^N \int_0^1\int_0^\infty\int_0^{2\pi}\Big[\phi({\bf v}+c_K(v_i, V_t^*({\bf v},{\bf w},\alpha),z,\varphi+
\varphi_{i,\alpha,t} ){\bf e}_i)-\phi({\bf v})\Big]  d\varphi dz d\alpha\\
=&\sum_{i=1}^N N^{-1}\sum_{j=1}^N \int_0^\infty\int_0^{2\pi}\Big[\phi({\bf v}+c_K(v_i, v_j, z,\varphi){\bf e}_i)-\phi({\bf v})\Big] d\varphi dz  \\
=&N^{-1}\sum_{i \ne j} \int_0^\infty\int_0^{2\pi}\Big[\phi({\bf v}+c_K(v_i, v_j, z,\varphi){\bf e}_i)-\phi({\bf v})\Big] d\varphi dz,
\end{align*}
because $c_K(v_i,v_i,z,\varphi)=0$. This is nothing but $\ca{L}_{N,K}\phi({\bf v})$, recall Lemma \ref{llxx}.

\vskip1mm

%Next, $\bb{E}\big[|V_t^{1}|^2\big]=m_2(f_0)$ follows from   a  standard computation:
%\begin{align*}
%\bb{E}\Big[|V_t^{1}|^2\Big]
%=&\bb{E}[|V_0^1|^2]+\frac{1}{N}\sum_{j=1}^N \int_0^t \int_0^\infty\int_0^{2\pi}
%\bb{E}\Big[|V_s^{1}+c_K(V_{s}^{1}, V_s^{j}, z,
%\varphi+\varphi_{1,\alpha,s})|^2-|V_s^{1}|^2\Big] d\varphi dz ds\\
%=&\bb{E}[|V_0^1|^2]+\frac{N-1}{N}\int_0^t \int_0^K
%\int_0^{2\pi} \bb{E}\Big[|c(V_{s}^{1}, V_s^{2}, z, \varphi)|^2
%+ 2V_s^{1}\cdot c(V_{s}^{1}, V_s^{2}, z, \varphi)\Big] d\varphi dz ds.
%\end{align*}
%We used the It\^{o} formula and that
%$\ca{L}_\alpha(V_t^*({\bf {V}}_{t},{\bf {W}}_{t}, \alpha))=N^{-1}\sum_{i=1}^N V_t^{i}$
%for the first equality;  the exchangeability, the $2\pi$-periodicity of $c_K=c\bbd{1}_{\{z\le K\}}$
%and that $c(V_{s}^{1}, V_s^{1}, z, \varphi)=0$ for the second one. Furthermore,
%Step 2 of the proof of Lemma \ref{estimategeneral} tells us that
%$$
%\bb{E}\Big[|V_t^{1}|^2\Big]
%=\bb{E}[|V_0^1|^2]+\frac{N-1}{N}\int_0^t \bb{E}\bigg[\Big( |V_{s}^{1}-V_s^{2}|^2
%- 2V_s^{1}\cdot (V_{s}^{1} - V_s^{2})\Big) \Phi_K(|V_{s}^{1}-V_s^{2}|)\bigg] ds.
%$$
%But the exchangeability implies that
%\begin{align*}
%\bb{E}\Big[|V_t^{1}|^2\Big]=\bb{E}[|V_0^1|^2]+\frac{N-1}{N}\int_0^t
%\bb{E}\bigg[\Big( |V_{s}^{1}-V_s^{2}|^2- V_s^{1}\cdot
%(V_{s}^{1} - V_s^{2}) - V_s^{2}\cdot (V_{s}^{2} - V_s^{1}) \Big)
% \Phi_K(|V_{s}^{1}-V_s^{2}|)\bigg] ds.
%\end{align*}
%Accordingly, $\bb{E}\Big[|V_t^{1}|^2\Big]=\bb{E}[|V_0^1|^2]$
%because the integrand is zero in the above formula.

Finally, we verify that $\sup_{[0,T]}\bb{E}\big[|V_t^{1}|^q\big]\le C_{T,q}$ if $f_0\in\ca{P}_q(\bb{R}^3)$
for some $q\ge2$: it immediately follows from  the It\^o formula, \eqref{q-increment} and exchangeability
that
\begin{align*}
\bb{E}\big[|V_t^{1}|^q\big]
&\le \bb{E}[|V_0^{1}|^q] + C_q \int_0^t\int_0^1\bb{E}\Big[1+|V_s^{1}|^q+|V_s^*({\bf {V}}_{s},  {\bf {W}}_{s}, \alpha)|^q\Big] \  d\alpha ds\\
&\le \bb{E}[|V_0^{1}|^q] + C_q N^{-1}\sum_{i=1}^{N} \int_0^t\bb{E}\Big[1+|V_s^{1}|^q+|V_s^{i}|^q\Big] \  ds\\
&\le \bb{E}[|V_0^{1}|^q] + C_q \int_0^t\bb{E}\Big[1+|V_s^{1}|^q\Big] \  ds,
\end{align*}
The Gr\"{o}nwall lemma allows us to complete the proof.
\epf

\brem
The exchangeability holds for the family $\{(W_t^i, V_t^{i})_{t\ge 0}, i=1,...,N\}$.
Indeed, the family $\{(W_t^i)_{t\ge 0}, i=1,...,N\}$ is i.i.d. by construction, so that the
exchangeability follows from the symmetry and pathwise uniqueness for \eqref{particle-system}.
\erem

\section{Bound of the Lebesgue norm of an empirical measure}\label{ttt}
An empirical measure cannot be in some $L^p$ space with $p>1$, so we will consider
a blob approximation, inspired by  \cite[Proposition 5.5]{Fournier:2015aa} and  \cite{MR3377068}.
But we deal with a jump process, so we need to overcome a few additional difficulties.

\vskip1mm

First, the following fact can be checked as 
Lemma 5.3 in \cite{Fournier:2015aa} (the norm and the step of the subdivision are different,
but this obviously changes nothing).

\blem\label{dicrete-norm}
Let $p\in(1,2)$ and $(f_t)_{t\geq 0} \in L^\infty \big([0,\infty),{\cal P}_2({\bb{R}^3})\big)\cap
L^1_{loc}\big([0,\infty), L^{p}({\bb{R}^3})\big)$
such that $m_2(f_t)=m_2(f_0)$ for all $t\geq 0$.
\benu[label=\emph{(\roman*)}]
\item There is a constant $M_p>0$, such that for all $t\ge0$, $\|f_t\|_{L^p}\ge M_p$.
\item For any $T>0$, we can find a subdivision $(t_{\ell}^N)_{\ell=0}^{K_N+1}$ satisfying $0=t_0^N<t_1^N<\dots<t_{K_N}^N \le T \le t_{K_N+1}^N$, such that $\sup_{\ell=0,...,K_N}(t_{\ell+1}^N-t_\ell^N)\le N^{-2}$ with $K_N \le 2TN^{2}$
and
\[ \int_0^Th_N(t)dt\le 2\int_0^T\|f_t\|_{L^p} dt, \]
with $h_N(t)=\sum_{\ell=1}^{K_N+1}\|f_{t_\ell^N}\|_{L^p}\bbd{1}_{\{t\in(t_{\ell-1}^N, t_\ell^N]\}}$.
\eenu
\elem

The goal of the section is to prove the following crucial fact.

\bprop\label{norm-bound}
Assume \eqref{con} for some $\gamma\in(-1,0)$, $\nu\in(0,1)$ with $\gamma+\nu>0$.
Let $q\geq 2$ such that $q>\gamma^2/(\gamma+\nu)$ and let $p\in(3/(3+\gamma),p_0(\gamma,\nu,q))\subset
(1,3/2)$. Consider $f_0 \in \ca{P}_q(\bb{R}^3)$ with a finite entropy and
$(f_t)_{t\ge0}\in L^\infty \big([0,\infty),{\cal P}_2({\bb{R}^3})\big)\cap
L^1_{loc}\big([0,\infty), L^{p}({\bb{R}^3})\big)$
the corresponding unique solution to \eqref{Bol} given by Theorem \emph{\ref{well-posedness}}.
Consider $(W_t^i)_{i=1,...,N, t\ge 0}$ the solution to \eqref{B-particle} and set
$\mu_{{\bf {W}}_t}^N=N^{-1}\sum_1^N \delta_{W^i_t}$.
Fix $\delta\in(0,1)$, set $\epsilon_N=N^{-(1-\delta)/3}$ and define
$\bar\mu_{{\bf {W}}_t}^N=\mu_{{\bf {W}}_t}^N*\psi_{\epsilon_N}$, where $\psi_\epsilon$ was defined in Notation
\emph{\ref{nnn}}. Finally, fix $T>0$ and consider $h_N$ built in Lemma \emph{\ref{dicrete-norm}}.
We have
 $$
 \bb{P}\Big(\forall t\in[0,T], \|\bar\mu_{{\bf {W}}_t}^N\|_{L^p}\le 13500 \big(1+h_N(t)\big)\Big)\ge 1- C_{T,q,\delta}N^{1-\delta q/3}.
 $$
\eprop

Throughout the section, we fix $N\ge 1$, $\delta\in(0,1)$, and $\epsilon_N=N^{-(1-\delta)/3}$ and
adopt the assumptions and notations of Proposition \ref{norm-bound}.  We  also put $r=p/(p-1)$.

\vskip1mm

In order to extend \cite[Proposition 5.5]{Fournier:2015aa}, it is necessary to study some
properties of the paths of the processes defined by \eqref{B-particle}.
Following Lemma 3.11 in \cite{Xu:2015aa}, we introduce, for each $i=1,\dots,N$,
\begin{align}\label{newprocess}
\widetilde{W}_t^i &=V_0^i+\int_0^t\int_0^1\int_0^\infty\int_0^{2\pi}c(W_{s-}^i,W_s^*(\alpha),z,\varphi)\bbd{1}_{\left\{|c(W_{s-}^i,W_s^*(\alpha),z,\varphi)|\le N^{-1/3}\right\}}M_i(ds,d\alpha,dz,d\varphi).
\end{align}

\blem\label{lem-small-jump}
 For all $T>0$,
$$
\bb{P}\Big[\sup_{[0,T]}|W_t^1|\le N^{\delta/3}, \sup_{s,t\in[0,T], |s-t|\le N^{-2}}|\widetilde{W}_t^1 - \widetilde{W}_s^1 |\ge \epsilon_N\Big]\le C_{T}N^2 e^{-N^{\delta/3}}.
$$
\elem
\bpf
Let us denote by $\tilde p$ the probability we want to bound.

\vskip1mm

{\it Step 1.}
We introduce
\begin{align*}
Z_t^1 &=\int_0^t\int_0^1\int_0^\infty\int_0^{2\pi} G\big(z/|W_{s-}^1-W_s^*(\alpha)|^\gamma\big)|W_{s-}^1-W_s^*(\alpha)| \\
& \hskip 5cm \times \bbd{1}_{\left\{G\big(z/|W_{s-}^1 - W_s^*(\alpha)|^\gamma\big)|W_{s-}^1 -W_s^*(\alpha)|/4\le N^{-1/3}\right\}} M_1(ds,d\alpha,dz,d\varphi).
\end{align*}
It is clear that $Z_t^1$ is almost surely increasing in $t$,  and  that a.s., for all $s, t\in[0, T]$,
\beq\label{pro-inequality}|\widetilde{W}_t^1 -\widetilde{W}_s^1 |\le |Z_t^1-Z_s^1|,
\eeq
since for any $v, v_*\in\bb{R}^3$ (recall \eqref{num})
$$
G\big(z/|v-v_*|^\gamma\big)|v-v_*|/4\le|c(v, v_* ,z,\varphi)|\le G\big(z/|v-v_*|^\gamma\big)|v-v_*|.$$
We now  consider the stopping time $\tau_N=\inf{\{t\ge0: |W_t^1|>N^{\delta/3}\}}$ and deduce from  \eqref{pro-inequality}  and the Markov inequality  that
\begin{align*}
\tilde p
&\le\bb{P}\Big[\sup_{[0,T]}|W_t^1|\le N^{\delta/3}, \sup_{s,t\in[0,T], |s-t|\le N^{-2}}|Z_t^1-Z_s^1|\ge \epsilon_N\Big]\\
&\le \bb{P}\Big[\sup_{s,t\in[0,T], |s-t|\le N^{-2}}|Z_{t\wedge\tau_N}^1-Z_{s\wedge\tau_N}^1|\ge \epsilon_N\Big].
\end{align*}
Since $[0, T] \subset \bigcup_{k=0}^{\lfloor N^{2}T\rfloor} [k/N^{2}, (k+1)/N^{2})$ and $Z_t^N$ is almost surely increasing in $t$, we deduce that on  $\{\sup_{s,t\in[0,T], |s-t|\le N^{-2}}|Z_{t\wedge\tau_N}^1-Z_{s\wedge\tau_N}^1|\ge \epsilon_N\}$, there exists $k\in\{0,1,...,\lfloor N^{2}T\rfloor \}$ for which  there holds $\big(Z_{((k+1)N^{-2})\wedge\tau_N}^1-Z_{(kN^{-2})\wedge\tau_N}^1\big)\ge \epsilon_N/3$. Hence,
\begin{align*}
\tilde p \le & \sum_{k=0}^{\lfloor N^{2}T\rfloor}\bb{P}\left[\Big(Z_{((k+1)N^{-2})\wedge\tau_N}^1-Z_{(kN^{-2})\wedge\tau_N}^1\Big)\ge
\frac{\epsilon_N}{3}\right]\\
&\le \sum_{k=0}^{\lfloor N^{2}T\rfloor}e^{-N^{\delta/3}}\bb{E}\left[\exp{\left\{3N^{1/3}\Big(Z_{((k+1)N^{-2})\wedge\tau_N}^1-Z_{(kN^{-2})\wedge\tau_N}^1\Big)\right\}}\right]\\
&=:\sum_{k=0}^{\lfloor N^{2}T\rfloor}e^{-N^{\delta/3}}I_k.
\end{align*}

{\it Step 2.} We now prove that $I_k$ is (uniformly) bounded, which will complete the proof.
We put
\[J_k(t)=:\bb{E}\left[\exp{\left\{3N^{1/3}\Big(Z_{(t+kN^{-2})\wedge\tau_N}^1-Z_{(kN^{-2})\wedge\tau_N}^1\Big)\right\}}\right].\]
It is obvious that $I_k=J_k(N^{-2})$. Applying  the  It\^{o} formula, we find
 \begin{align*}
 J_k(t)
 &=1+2\pi \bb{E}\bigg[\int_{(kN^{-2})\wedge\tau_N}^{(t+kN^{-2})\wedge\tau_N} \int_0^1 \int_0^\infty  \exp{\left\{3N^{1/3}\Big(Z_s^1-Z_{(kN^{-2})\wedge\tau_N}^1\Big)\right\}} \\
 &\quad \times \Big(e^{3N^{1/3}G\big(z/|W_{s}^1-W_s^*(\alpha)|^\gamma\big)|W_{s}^1-W_s^*(\alpha)|}-1\Big) \bbd{1}_{\left\{G\big(z/|W_{s}^1 - W_s^*(\alpha)|^\gamma\big)|W_{s}^1 -W_s^*(\alpha)|/4\le N^{-1/3}\right\}} dz d\alpha ds\bigg].
 \end{align*}
 Since  $3N^{1/3}G\big(z/|W_{s}^1-W_s^*(\alpha)|^\gamma\big)|W_{s}^1-W_s^*(\alpha)|\le 12$
(thanks to the indicator function), we have
$$
 e^{3N^{1/3}G\big(z/|W_{s}^1-W_s^*(\alpha)|^\gamma\big)|W_{s}^1-W_s^*(\alpha)|}-1 \le CN^{1/3}G\big(z/|W_{s}^1-W_s^*(\alpha)|^\gamma\big)|W_{s}^1-W_s^*(\alpha)|
 $$ for a positive constant $C$. Then using \eqref{con1}, we see that
 $$
 \bbd{1}_{\left\{G\big(z/|W_{s}^1 - W_s^*(\alpha)|^\gamma\big)|W_{s}^1 -W_s^*(\alpha)|/4\le N^{-1/3}\right\}} \le \bbd{1}_{\left\{z \ge CN^{\nu/3}|W_{s}^1-W_s^*(\alpha)|^{\gamma+\nu}-|W_{s}^1-W_s^*(\alpha)|^{\gamma}\right\}}.
 $$
Hence,
\begin{align*}
 J_k(t) &\le  1+ CN^{1/3}\bb{E}\bigg[\int_{(kN^{-2})\wedge\tau_N}^{(t+kN^{-2})\wedge\tau_N} \int_0^1 \int_0^\infty  \exp{\left\{3N^{1/3}\Big(Z_s^1-Z_{(kN^{-2})\wedge\tau_N}^1\Big)\right\}} \\
 & \quad \times \Big(1+z/|W_{s}^1 - W_s^*(\alpha)|^\gamma\Big)^{-1/\nu} |W_{s}^1 -W_s^*(\alpha)| \bbd{1}_{\left\{z \ge CN^{\nu/3}|W_{s}^1-W_s^*(\alpha)|^{\gamma+\nu}-|W_{s}^1-W_s^*(\alpha)|^{\gamma}\right\}} dz d\alpha ds\bigg].
\end{align*}
But, we have
\begin{align*}
&|W_{s}^1-W_s^*(\alpha)| \int_0^\infty \Big(1+z/|W_{s}^1 - W_s^*(\alpha)|^\gamma\Big)^{-1/\nu} \bbd{1}_{\left\{z \ge CN^{\nu/3}|W_{s}^1-W_s^*(\alpha)|^{\gamma+\nu}-|W_{s}^1-W_s^*(\alpha)|^{\gamma}\right\}} dz\\
=& C N^{-(1-\nu)/3}|W_{s}^1-W_s^*(\alpha)|^{\nu+\gamma} \\
\le& CN^{-(1-\nu)/3}(1+|W_{s}^1|^2+|W_s^*(\alpha)|^2 )
\end{align*}
since $\gamma+\nu\in(0,1)$.
Using now that $\int_0^1|W_s^*(\alpha)|^2 d\alpha=m_2(f_0)$ and that $|W_s^1|\le N^{\delta/3}$ for all $s\le \tau_N$,
we conclude that
 \begin{align*}
 J_k(t)&\le 1+CN^{\nu/3}(1+m_2(f_0)+N^{2\delta/3})\int_0^t J_k(s)ds\\
 &\le 1+CN^{(\nu+2\delta)/3}\int_0^t J_k(s)ds.
 \end{align*}
It follows from the Gr\"{o}nwall lemma that $J_k(t)\le \exp{(CN^{(\nu+2\delta)/3}t)}$, and thus
that $I_k=J_k(N^{-2})$ is uniformly bounded, because $(\nu+2\delta)/3< 2$
(recall that $\nu\in(0,1)$ and $\delta\in(0,1)$).
\epf

Next, we study the {\it large} jumps of $(W_t^1)_{t\geq 0}$.

\blem\label{lem-big-jump}
There exists $C>0$ such that for any $\ell\in\{1,..., K_N+1\}$,
$$
\bb{P}\Big[\exists~ t\in(t_{\ell-1}^N, t_{\ell}^N]: |\Delta W^1_t|>N^{-1/3}\Big] \le CN^{-(2-\nu/3)}.
$$
\elem

\bpf
Let us fix $\ell$ and set $A=\{\exists~ t\in(t_{\ell-1}^N, t_{\ell}^N]: |\Delta W^1_t|>N^{-1/3}\}$.
After noting that
\begin{align*}
A&=\left\{\int_{t_{\ell-1}^N}^{t_\ell^N}\int_0^1\int_0^\infty\int_0^{2\pi}\bbd{1}_{\left\{|c(W_{s-}^i,W_s^*(\alpha),z,\varphi)|> N^{-1/3}\right\}}M_1(ds, d\alpha, dz, d\varphi)\ge 1\right\},
\end{align*}
we have
\begin{align*}
\bb{P}(A)
&\le \bb{E}\left[\int_{t_{\ell-1}^N}^{t_\ell^N}\int_0^1\int_0^\infty\int_0^{2\pi}\bbd{1}_{\left\{|c(W_{s-}^1,W_s^*(\alpha),z,\varphi)|> N^{-1/3}\right\}}M_1(ds, d\alpha, dz, d\varphi) \right]
\end{align*}
by  the Markov inequality.
But,  \eqref{num} and \eqref{con1} tell us  that $|c(v,v_*,z, \varphi)|
\leq C(1+z/|v-v_*|^\gamma)^{-1/\nu}|v-v_*|$.  Hence,
 \begin{align*}
\bb{P}(A)
 & \le  2\pi \bb{E}\Big[\int_{t_{\ell-1}^N}^{t_\ell^N}\int_0^1\int_0^\infty \bbd{1}_{\left\{C ( 1+z/|W_s^1 -W_s^*(\alpha)|^\gamma)^{-1/\nu} |W_s^1-W_s^*(\alpha)|> N^{-1/3}\right\}} dz d\alpha  ds\Big]\\
 & \le 2\pi \bb{E}\Big[\int_{t_{\ell-1}^N}^{t_\ell^N}\int_0^1\int_0^\infty \bbd{1}_{\left\{ z< CN^{\nu/3} |W_s^1-W_s^*(\alpha)|^{\gamma+\nu} \right\}}  dz d\alpha  ds\Big]\\
 &= C N^{\nu/3} \bb{E}\Big[\int_{t_{\ell-1}^N}^{t_\ell^N}\int_0^1|W_s^1-W_s^*(\alpha)|^{\gamma+\nu} d\alpha  ds\Big].
 \end{align*}
Finally, using that  $|W_s^1-W_s^*(\alpha)|^{\gamma+\nu}\le 1+|W_s^1|^2+|W_s^*(\alpha)|^2$
and that $\int_0^1|W_s^*(\alpha)|^2 d\alpha=\bb{E}[|W_s^1|^2]<\infty$,
we conclude that $\bb{P}(A)\leq CN^{\nu/3}(t_{\ell+1}^N-t_\ell^N)\le C N^{\nu/3-2}$ as desired.
\epf

\blem\label{event}
For $\ell=1, ..., K_N+1$, we introduce
\beq\label{Number}
I_{\ell}=\{i\in\{1,...,N\}: \exists~ t\in (t_{\ell-1}^N, t_{\ell}^N] ~\text{such that}~ |\Delta W_t^i|>N^{-1/3}\},
\eeq
and the event
\begin{align*}
\Omega_{T,N}^1= & \left\{\forall i\in\{1,..,N\},\; \sup_{[0,T]}|W_t^i|\le N^{\delta/3}
\text{ and }
\sup_{s,t\in[0,T], |s-t|\le N^{-2}}|\widetilde{W}_t^i - \widetilde{W}_s^i|\le \epsilon_N\right\}\\
&\ \bigcap \, \left\{\forall \ell=1,...,K_N+1, \;\#(I_\ell) \le N \epsilon_N^{3/r}\right\}.
\end{align*}
Then we have
\[\bb{P}[\Omega_{T,N}^1] \ge 1-C_{T,q,\delta}N^{1-q\delta/3}.\]
\elem

\bpf
We write $\Omega_{T,N}^1=\Omega_{T,N}^{1,1} \cap \Omega_{T,N}^{1,2}$, where
\begin{align*}
\Omega_{T,N}^{1,1} &:= \left\{\forall i\in\{1,..,N\},\; \sup_{[0,T]}|W_t^i|\le N^{\delta/3}\text{ and } \sup_{s,t\in[0,T], |s-t|\le N^{-2}}|\widetilde{W}_t^i - \widetilde{W}_s^i|\le \epsilon_N\right\},\\
\Omega_{T,N}^{1,2} &:= \left\{\forall \ell=1,...,K_N+1,\;  \#(I_\ell) \le N \epsilon_N^{3/r}\right\},
\end{align*}

{\it Step 1.} Here we estimate $\bb{P}[(\Omega_{T,N}^{1,1})^c]$.
Using  the Markov inequality,  \eqref{lem-momen}   and   Lemma \ref{lem-small-jump}, we get
\begin{align*}
\bb{P}[(\Omega_{T,N}^{1,1})^c]
&\le N \, \bb{P}\Big[\Big\{\sup_{[0,T]}|W_t^{1}|\le N^{\delta/3} \text{ and } \sup_{|s-t|\le N^{-2}} |\widetilde{W}_t^{1}-\widetilde{W}_s^1|\le \epsilon_N\Big\}^c\Big]\\
&= N\, \bb{P}\Big[\sup_{[0,T]}|W_t^{1}|\ge N^{\delta/3}\Big] + N \,\bb{P}\Big[\sup_{[0,T]}|W_t^{1}|\le N^{\delta/3} \text{ and } \sup_{|s-t|\le N^{-2}} |\widetilde{W}_t^{1}-\widetilde{W}_s^1|\ge \epsilon_N\Big]\\
&\le N\,  \bb{E}\Big[\sup_{[0,T]} |W_t^1|^q \Big] \, N^{-q\delta/3} + C_{T} \, N^3 e^{-N^{\delta/3}}\le C_{T,q} \,N^{1-q\delta/3}.
\end{align*}

\vskip1mm 

{\it Step 2.} We now prove that  $\bb{P}[(\Omega_{T,N}^{1,2})^c]\le C_T\exp{(-N^{\delta})}$. For any fixed $\ell\in\{1,...,K_N+1\}$, we introduce  $A_N^\ell=\{\exists~ t\in(t_{\ell-1}^N, t_{\ell}^N]: |\Delta W^1_t|>N^{-1/3}\}$.  Then we observe that $\#(I_\ell)$ follows a Binomial distribution with parameters $N$ and $\bb{P}(A_N^\ell)$. Using again  the Markov inequality, we observe that
\begin{align*}
\bb{P}[(\Omega_{T,N}^{1,2})^c]\le  \sum_{\ell=1}^{K_N+1}\bb{P}\Big[  \#(I_\ell) \ge N\epsilon_N^{3/r}\Big]\le \sum_{\ell=1}^{K_N+1}\bb{E}[\exp \big( \#(I_\ell) \big)] \exp{(-N\epsilon_N^{3/r})}.
\end{align*}
But, \[\bb{E}[\exp{( \#(I_\ell))}]=\exp{\big(N \log(1  + (e-1) \bb{P}(A_N^{\ell}))\big)}\le \exp{\big(N(e-1) \bb{P}(A_N^{\ell})\big)}.
\]
Hence,
\[\bb{P}[(\Omega_{T,N}^{1,2})^c]\le  \sum_{\ell=1}^{K_N+1}\exp{\big(N(e-1) \bb{P}(A_N^{\ell})\big)}\exp{(-N\epsilon_N^{3/r})}.\]
We know  from Lemma \ref{lem-big-jump} that
$\bb{P}(A_N^\ell)\le C N^{-(2-\nu/3)}$, hence $N\bb{P}(A_N^\ell) \le CN^{-1+\nu/3}\le C$. We thus  deduce  that
$$
\bb{P}[(\Omega_{T,N}^{1,2})^c]\le C (K_N+1)\exp(-N\epsilon_N^{3/r}) \le C(2TN^2+1)\exp(-N\e_N^{3/r})\le C_T \exp(-N^{\delta}),
$$
since $N\e_N^{3/r}=N^{1/p+\delta/r}$ and since $1/p+\delta/r>\delta$. This ends the proof.
\epf
We now give the
\bpf[Proof of Proposition \ref{norm-bound}]
Consider the partition $\scr{P}_N$ of  $\bb{R}^3$ in cubes with side length $\e_N$ and its subset $\scr{P}_N^\delta$  consisting of cubes that have non-empty intersection with  $B(0, N^{\delta/3})$. Then we observe that $\# (\scr{P}_N^\delta)\le (2(N^{\delta/3}+\epsilon_N)\epsilon_N^{-1})^3\le 64 N^\delta\epsilon_N^{-3}=64N$. We split the proof into several steps.

\vskip1mm 

{\it Step 1.}   For  $(x_1,..., x_N)\in(B(0,N^{\delta/3}))^N$ and $(y_1,..., y_N)\in(B(0,N^{\delta/3}))^N$, we set
\[ I= \{i\in\{1,.., N\}: |x_i-y_i|> \epsilon_N\},\]
and  denote the empirical measure of  ${\bf y}=(y_1,..., y_N)\in(\bb{R}^3)^N$ by $\mu_{{\bf y}}^N=N^{-1}\sum_{i=1}^N\delta_{y_i}$.
The goal of this step is to show that
\[\| \mu_{{\bf y}}^N*\psi_{\epsilon_N} \|_{L^p}\le \big(\frac{3}{4\pi}\big)^{1/r} \frac{\#(I)}{N\epsilon_N^{3/r}} + 3375 \Big(N^{-p}\epsilon_N^{-3(p-1)}\sum_{D\in\scr{P}_N^\delta}(\#\{i\in\{1,...,N\} : x_i\in D\})^p\Big)^{1/p}.\]
Indeed, recalling that $\psi_\epsilon(x)=(3/(4\pi\epsilon^3))\bbd{1}_{\{|x|\le \epsilon\}}$, we observe that
\begin{align*}
\mu_{{\bf y}}^N*\psi_{\epsilon_N}(v)
& = \frac{1}{N} \sum_{i=1}^{N} \psi_{\epsilon_N}(v-y_i)\bbd{1}_{\{|x_i-y_i|>\epsilon_N\}}+N^{-1} \sum_{i=1}^{N} \psi_{\epsilon_N}(v-y_i)\bbd{1}_{\{|x_i-y_i| \le \epsilon_N\}}\\
& =  \frac{1}{N} \sum_{i\in I}\psi_{\epsilon_N}(v-y_i)
+ \frac{3}{4\pi N\epsilon_N^3} \# \big\{ i\in\{1,...,N\}: y_i\in B(v,\epsilon_N), |y_i-x_i|\le \epsilon_N \big\} \\
& \le \frac{1}{N} \sum_{i\in I} \psi_{\epsilon_N}(v-y_i) + \frac{3}{4\pi N\epsilon_N^3}\# \big \{ i\in\{1,...,N\}: x_i\in B(v , 2\epsilon_N)\big \}.
\end{align*}
Hence,
$$
\mu_{{\bf y}}^N*\psi_{\epsilon_N}(v) \le  \frac{1}{N} \sum_{i\in I} \psi_{\epsilon_N}(v-y_i) +\frac{3}{4\pi N\epsilon_N^3} \sum_{D\in\scr{P}_N^\delta}\# \big \{  i\in\{1,...,N\} :  x_i\in D \big \} \bbd{1}_{\{D\cap B(v, 2\epsilon_N) \neq \emptyset\}}.
$$
We then deduce that
\begin{align*}
&\|\mu_{{\bf y}}^N*\psi_{\epsilon_N}\|_{L^p}\\
&\le \frac{1}{N} \| \sum_{i\in I} \psi_{\epsilon_N}(\cdot -y_i)\|_{L^p}
+ \frac{3}{4\pi N\epsilon_N^3} \| \sum_{D\in\scr{P}_N^\delta}\# \{ i\in\{1,...,N\}:  x_i\in D\} \bbd{1}_{\{D\cap B(\cdot, 2\epsilon_N)\neq \emptyset\}} \|_{L^p}.
 \end{align*}
Since $\| \psi_{\epsilon_N}(\cdot -y_i) \|_{L^p}=(\frac{3}{4\pi})^{1/r} \epsilon_N^{-3/r}$, we have
\begin{align*}
\frac{1}{N} \| \sum_{i\in I} \psi_{\epsilon_N}(\cdot -y_i) \|_{L^p}
 \le \frac{1}{N}\sum_{i\in I} \| \psi_{\epsilon_N}(\cdot -y_i) \|_{L^p} &\le \big(\frac{3}{4\pi} \big)^{1/r} \frac{\#(I)}{N\epsilon_N^{3/r}}.
\end{align*}
On the other hand,  let $A := \| \sum_{D\in\scr{P}_N^\delta}\# \{i\in\{1,...,N\} :  x_i\in D\} \bbd{1}_{\{D\cap B(\cdot, 2\epsilon_N)\neq \emptyset\}} \|_{L^p}$, then
\begin{align*}
A^p
&=\int_{\bb{R}^3}\Big(\sum_{D\in\scr{P}_N^\delta} \# \{i:  x_i\in D\} \bbd{1}_{\{D\cap B(v, 2\epsilon_N)\neq \emptyset\}}\Big)^p dv\\
&=\int_{\bb{R}^3}\Big(\sum_{D,D^\prime\in\scr{P}_N^\delta} \# \{i:  x_i\in D\}\# \{i:  x_i\in D^\prime\} \bbd{1}_{\{D\cap B(v, 2\epsilon_N)\neq \emptyset, D^\prime\cap B(v, 2\epsilon_N)\neq \emptyset\}}\Big)^{p/2} dv\\
&\le \int_{\bb{R}^3} \sum_{D,D^\prime\in\scr{P}_N^\delta} \big(\# \{i:  x_i\in D\} \big)^{p/2} \big(\# \{i:  x_i\in D^\prime\} \big)^{p/2} \bbd{1}_{\{D\cap B(v, 2\epsilon_N)\neq \emptyset, D^\prime\cap B(v, 2\epsilon_N)\neq \emptyset\}} dv
\end{align*}
because $p\in(1,2)$.
From  $x^2+y^2\ge 2xy$ and a symmetry argument, we see that
\begin{align*}
A^p \le \sum_{D\in\scr{P}_N^\delta} (\# \{i: x_i\in D\})^p \int_{\bb{R}^3}
 \bbd{1}_{\{D\cap B(v, 2\epsilon_N)\neq \emptyset\}} \sum_{D^\prime\in\scr{P}_N^\delta}  \bbd{1}_{\{D^\prime \cap B(v, 2\epsilon_N)\neq \emptyset \}} dv.
\end{align*}
But,  for each $v\in\bb{R}^3$,
$\sum_{D^\prime\in\scr{P}_N^\delta}  \bbd{1}_{\{D^\prime \cap B(v, 2\epsilon_N)\neq \emptyset \}}=\#\{D^\prime\in\scr{P}_N^\delta: D^\prime \cap B(v, 2\epsilon_N)\neq \emptyset\} \le 5^3.$
And for each $D\in\scr{P}_N^\delta$, $\{v\in\bb{R}^3: D\cap B(v, 2\epsilon_N)\neq \emptyset\}$ is included by a ball of radius $3\e_N$. Therefore,
 $\int_{\bb{R}^3}
 \bbd{1}_{\{D\cap B(v, 2\epsilon_N)\neq \emptyset\}} dv \le 4\pi(3\epsilon_N)^3/3$. Hence,
 $$
 A^p \le \frac{5^3 4\pi (3\e_N)^3}{3}\sum_{D\in\scr{P}_N^\delta} \big(\# \{i: x_i\in D\} \big)^p.
 $$
 Consequently,
 \begin{align*}
 \|\mu_{{\bf y}}^N*\psi_{\epsilon_N}(v)\|_{L^p}
 &\le \big(\frac{3}{4\pi} \big)^{1/r} \frac{\#(I)}{N\epsilon_N^{3/r}}
 + \frac{3}{4\pi N\epsilon_N^3} A \\
 &\le \big(\frac{3}{4\pi} \big)^{1/r} \frac{\#(I)}{N\epsilon_N^{3/r}}  + \big(\frac{3}{4\pi} \big)^{1/r} (15)^{3/p}
\Big(N^{-p}\epsilon_N^{-3(p-1)}\sum_{D\in\scr{P}_N^\delta} \big(\#\{ i: x_i\in D\} \big)^p\Big)^{1/p}.
 \end{align*}
Since $(15)^{3/p}\leq 15^3=3375$, this ends the step.

\vskip1mm 

{\it Step 2.} In this step, we extend  the  proof of  \cite[Step 3-Proposition 5.5]{Fournier:2015aa}
to show that  there are some constants $C>0$ and $c>0$ (depending on $\delta$ and $M_p$ of
Lemma \ref{dicrete-norm}) such that for all  fixed $t\in[0, T+1]$,
 \[\bb{P}[(\Omega_{t,N}^2)^c]\le C \exp{(-cN^{\delta/r})},\]
 where \[\Omega_{t,N}^2=\left\{N^{-p}\epsilon_N^{-3(p-1)}\sum_{D\in\scr{P}_N^\delta}\Big(\#\{i\in\{1,...,N\} : W_t^i\in D\}\Big)^p \le 2^{p+1} \| f_t \|_{L^p}^p \right\}.\]
 To this end, we introduce, for $D\in\scr{P}_N^\delta$, $A_D=\#\{i: W_t^i\in D\}$. Then $A_D\sim B(N, f_t(D))$
and we have
 \beq\label{Prop-proba}
 \bb{P}(A_D\ge x)\le \exp(-x/8) \quad \text{for\, all}\quad  x\ge 2Nf_t(D).
 \eeq
Indeed, $\bb{P}(A_D\ge x)\le e^{-x} \bb{E}[\exp(A_D)]=e^{-x}\exp[N\log(1+f_t(D)(e-1))] \le e^{-x}\exp[N(e-1)f_t(D)]$. If  $x\ge 2Nf_t(D)$, we thus have
$$
\bb{P}(A_D\ge x) \le \exp[-x+x(e-1)/2] \le \exp(-x/8).
$$
Next, it follows  from  the H\"{o}lder inequality that
$$
\|f_t\|_{L^p}^p \ge \sum_{D\in\scr{P}_N^\delta} \int_D |f_t(v)|^p dv \ge  \e_N^{-3p/r} \sum_{D\in\scr{P}_N^\delta}(f_t(D))^p.
$$
On the other hand,  we observe from $\# (\scr{P}_N^\delta)\le 64 N^\delta\epsilon_N^{-3}$ that
$$
\|f_t\|_{L^p}^p  \ge  64^{-1}N^{-\delta}\e_N^3 \sum_{D\in\scr{P}_N^\delta} \|f_t\|_{L^p}^p.
$$
Using the two previous inequalities, we find that
$$
2^{p+1}\|f_t\|_{L^p}^p  \ge \sum_{D\in\scr{P}_N^\delta} \big(2^p\e_N^{-3p/r}(f_t(D))^p+
2^p 64^{-1}N^{-\delta}\e_N^3\|f_t\|_{L^p}^p \big).
$$
Consequently, on  $(\Omega_{t,N}^2)^c$, we have
$$
\sum_{D\in\scr{P}_N^\delta} A_D^p >N^p\e_N^{3(p-1)}2^{p+1}\|f_t\|_{L^p}^p \geq N^p\e_N^{3(p-1)}
\sum_{D\in\scr{P}_N^\delta} \big(2^p\e_N^{-3p/r}(f_t(D))^p+2^p 64^{-1}N^{-\delta}\e_N^3\|f_t\|_{L^p}^p \big),
$$
so that there is
at least one $D\in\scr{P}_N^\delta$ with $A_D^p \ge N^p \e_N^{3(p-1)}
\big[2^p\e_N^{-3p/r}(f_t(D))^p+2^p 64^{-1}N^{-\delta}\e_N^3\|f_t\|_{L^p}^p \big] $. Hence,
 $$
\bb{P}[(\Omega_{t,N}^2)^c] \le \sum_{D\in\scr{P}_N^\delta} \bb{P}\Big(A_D \ge N \e_N^{3/r}
\big[2^p\e_N^{-3p/r}(f_t(D))^p+2^p 64^{-1}N^{-\delta}\e_N^3\|f_t\|_{L^p}^p\big]^{1/p}\Big).
 $$
 But we can apply \eqref{Prop-proba}, because $x_N:=N \e_N^{3/r} \big[2^p\e_N^{-3p/r}(f_t(D))^p+2^p64^{-1}N^{-\delta}\e_N^3\|f_t\|_{L^p}^p\big]^{1/p}$ enjoys
the property that $x_N\ge N \e_N^{3/r} [2^p\e_N^{-3p/r}(f_t(D))^p]^{1/p}=2Nf_t(D)$:
  $$ \bb{P}[(\Omega_{t,N}^2)^c] \le \sum_{D\in\scr{P}_N^\delta} \exp(-x_N/8).$$
Using that  $x_N \ge N \e_N^{3/r}(2^p 64^{-1}N^{-\delta}\e_N^3\|f_t\|_{L^p}^p)^{1/p}= cN^{\delta/r}\|f_t\|_{L^p}$,
that $\#(\scr{P}_N^\delta) \le 64N$ and that $\|f_t\|_{L^p}\ge M_p$, we deduce that
  $$
  \bb{P}[(\Omega_{t,N}^2)^c] \le \sum_{D\in\scr{P}_N^\delta} \exp(- cN^{\delta/r}\|f_t\|_{L^p}/8)
\le 64N \exp(-c M_pN^{\delta/r}/8) \le C\exp(-cM_p N^{\delta/r}/10).
  $$
  This ends the step.

  \vskip1mm 
  
  {\it Step 3.}  We finally consider the event
 \[\Omega_{T,N}=\Omega_{T,N}^1\cap(\cap_{\ell=1}^{K_N+1}\Omega_{t_\ell^N, N}^2),
 \]
 where $\Omega_{T,N}^1$ is defined in Lemma \ref{event}, and the sequence $(t_{\ell}^N)_{\ell=0}^{K_N+1}$ satisfying  $0=t_0^N < t_1^N<...<t_{K_N}^N \le T\le T_{K_N+1}^N$, with $K_N\le 2TN^2$ and $\sup_{i=0,...,K_N}(t_{\ell+1}^N-t_\ell^N)\le N^{-2}$  is built  in Lemma \ref{dicrete-norm}. We also recall that $h_N(t)=\sum_{\ell=1}^{K_N+1}\|f_{t_\ell^N}\|_{L^p}\bbd{1}_{\{t\in(t_{\ell-1}^N, t_\ell^N]\}}$.

\vskip1mm

According to Lemma \ref{event} and Step 2, we see that
\[\bb{P}[\Omega_{T,N}^c] \leq  \bb{P}[(\Omega_{T,N}^1)^c]+\sum_{\ell=1}^{K_N+1}\bb{P}[(\Omega_{t_\ell^N, N}^2)^c]\le C_{T,q,\delta}N^{1-q\delta/3}+C(K_N+1)\exp{(-cN^{\delta/r})}\le C_{T,q,\delta}N^{1-q\delta/3}.\]

\vskip1mm

 Finally, we show that  on $\Omega_{T,N}$, for all $t\in[0,T]$, $\|\bar\mu_{{\bf {W}}_t}^N\|_{L^p}\le 13500 (1+h_N(t))$.
Recall that $\widetilde{W}_t^i$ is defined by \eqref{newprocess} and that $I_\ell$ is given by  \eqref{Number}, we have
 \benu[label=(\roman*)]
 \item for all $i=1,...,N$, and for all $t\in [0, T+1]$, $W_t^i\in B(0, N^{\delta/3})$
(according to $\Omega_{T,N}^1$);
 \item for all $\ell=1,..., K_N+1$, all $t\in(t_{\ell-1}^N, t_{\ell}^N]$, all $i\in\{1,...,N\}\setminus I_\ell$,
$|W_t^i-W_{t_\ell^N}^i|=|\widetilde{W}_t^i-\widetilde{W}_{t_\ell^N}^i|\le \epsilon_N$, and
$ \#(I_\ell) \le N \epsilon_N^{3/r}$ (by definition of $\widetilde W^i$ and $I_\ell$ and thanks to $\Omega_{T,N}^1$);
 \item For all $\ell=1,..., K_N+1$, $N^{-p}\epsilon_N^{-3(p-1)}\sum_{D\in\scr{P}_N^\delta}\Big(\#\big\{i\in\{1,...,N\} : \, W_{t_\ell^N}^i\in D\big \}\Big)^p\le 2^{p+1} \| f_{t_\ell^N} \|_{L^p}^p$ (according to $\cap_{\ell=1}^{K_N+1}\Omega_{t_\ell^N, N}^2$).
 \eenu

Using Step 1 with  $\bar\mu_{{\bf {W}}_t}^N=\mu_{{\bf {W}}_t}^N*\psi_{\epsilon_N}$,
we deduce that on $\Omega_{T,N}$, for all $t\in[0,T]$, choosing $\ell$ such that $t\in(t_{\ell-1}^N, t_\ell^N]$, we
have
\begin{align*}
\|\bar\mu_{{\bf {W}}_t}^N\|_{L^p}\le&  \big(\frac{3}{4\pi}\big)^{1/r} \frac{\#(I_\ell)}{N\epsilon_N^{3/r}}
+ 3375 \Big(N^{-p}\epsilon_N^{-3(p-1)}\sum_{D\in\scr{P}_N^\delta}(\#\{i\in\{1,...,N\} : W^i_{t^N_{\ell}}\in D\})^p\Big)^{1/p}\\
\leq & 1+ 3375. 2^{(p+1)/p}\| f_{t_\ell^N} \|_{L^p}\\
=& 1+ 3375. 2^{(p+1)/p}h_N(t).
\end{align*}
This completes the proof, since  $3375. 2^{(p+1)/p}\leq 3375.4=13500$.
\epf

\section{Estimate of the Wasserstein distance}

This last section is devoted to the proof of Theorem \ref{main-result}.
In the whole section, we assume \eqref{con} for some $\gamma\in(-1,0)$, $\nu\in(0,1)$ with $\gamma+\nu>0$.
We consider $q>6$ such that  $q>\gamma^2/(\gamma+\nu)$, $f_0\in\ca{P}_q({\bb{R}^3})$
with a finite entropy, and $(f_t)_{t\ge0}$ the  unique weak solution to  \eqref{Bol}
given by Theorem \ref{well-posedness}. We fix  $p\in (3/(3+\gamma),p_0(\gamma,\nu,q))$
and  know that $(f_t)_{t\geq 0}\in L^\infty\big([0,\infty), \ca{P}_2({\bb{R}^3})\big)\cap
L^1_{loc}\big([0,\infty), L^{p}({\bb{R}^3})\big)$.

\vskip1mm

We fix $N\geq 1$, $K\geq 1$ and put $\e_N=N^{-(1-\delta)/3}$ with $\delta=6/q$.
Consider $(V_t^{i})_{t\ge0}$  for $i=1,\dots,N$, defined by \eqref{particle-system}
with the choice $\e=\e_N$.
We know by Lemma \ref{aaaa} that  $(V_t^{i})_{i=1,\dots,N,t\ge0}$ is a Markov process with generator
${\cal L}_{N,K}$, see \eqref{generator-cut}, starting from $(V_0^{i})_{i=1,\dots,N}$, which is an i.i.d. family
of $f_0$-distributed random variables. We set $\mu^N_{{\bf{V}}_t}=N^{-1}\sum_{1}^N \delta_{V^i_t}$.
So the goal of the section is to prove that
\begin{align}\label{ggooal}
\sup_{[0, T]}\bb{E}[\ca{W}_2^2(\mu^N_{{\bf{V}}_t}, f_t)]\le C_{T,q}\Big(N^{-(1-6/q)(2+2\gamma)/3}+K^{1-2/\nu}+N^{-1/2}\Big).
\end{align}
We consider $(W_t^i)_{t\ge0}$, for $i=1,\dots,N$ defined by \eqref{B-particle} and recall that
for all $t\geq0$, the family $(W_t^i)_{i=1,\dots, N}$ is i.i.d. and $f_t$-distributed.

\vskip1mm

First, we introduce the following shortened notations:
\begin{align*}
&c_W(s):=c(W_{s}^{1},W_s^*(\alpha), z,\varphi),\\
&c_W^{N}(s):=c(W_{s}^{1},W_s^{*,\epsilon_N}(\alpha), z,\varphi+\varphi_{1,\alpha, s}^1),\\
&c_{V}^{N}(s):=c(V_{s}^{1},V_s^{*,\epsilon_N}({\bf {V}}_{s},  {\bf {W}}_{s},  \alpha), z,\varphi+\varphi_{1,\alpha, s}^1+\varphi_{1,\alpha, s}^2),\\
&c_{K,V}^{N}(s):=c_K(V_{s}^{1},V_s^{*,\epsilon_N}({\bf {V}}_{s},  {\bf {W}}_{s},  \alpha), z,\varphi+\varphi_{1,\alpha, s}^1+\varphi_{1,\alpha, s}^2),\\
&c_{K,V}(s):=c_K(V_{s}^{1}, V_s^*({\bf {V}}_{s},  {\bf {W}}_{s}, \alpha), z, \varphi+\varphi_{1,\alpha,s}),
\end{align*}
with the notations of
Section 4.
Let us now prove an intermediate result.

\blem\label{intermediate-lemma}
There is $C>0$ such that a.s.,
\begin{align*}
&I_0^N(s)+I_1^N(s)+I_2^N(s)+I_3^N(s) \\
\le& C \e_N^{2+2\gamma} + C|W_s^1-V_s^{1}|^2+ C
K^{1-2/\nu}\int_0^1 |W_s^1-W_s^{*,\epsilon_N}(\alpha)|^{2+2\gamma/\nu} d\alpha\\
&+C\int_0^1\Big(|W_s^1-V_s^{1}|^2+|W_s^{*}(\alpha)-V_s^{*}({\bf {V}}_{s},  {\bf {W}}_{s},  \alpha)|^2\Big)|W_s^1-W_s^{*,\epsilon_N}(\alpha)|^\gamma d\alpha.
\end{align*}
where
\begin{align*}
& I_0^N(s):=\int_0^1\int_0^\infty\int_0^{2\pi}\Big(2(W_{s}^1-V_{s}^{1})\cdot(c_W^{N}(s)-c_{K,V}^{N}(s))+|c_W^N(s)-c_{K,V}^N(s)|^2\Big)d\varphi dzd\alpha,\\
& I_1^N(s):=\int_0^1\int_0^\infty\int_0^{2\pi}2(W_{s}^1-V_{s}^{1})\cdot\big(c_W(s)-c_W^{N}(s)+c_{K,V}^{N}(s)-c_{K,V}(s)\big) d\varphi dzd\alpha,\\
& I_2^N(s):=\int_0^1\int_0^\infty\int_0^{2\pi}|c_W(s)-c_W^{N}(s)+c_{K,V}^{N}(s)-c_{K,V}(s)|^2 d\varphi dzd\alpha,\\
& I_3^N(s):=\int_0^1\int_0^\infty\int_0^{2\pi} 2\big(c_W^N(s)-c_{K,V}^N(s)\big)\cdot\big(c_W(s)-c_W^{N}(s)+c_{K,V}^{N}(s)-c_{K,V}(s)\big) d\varphi dzd\alpha.
\end{align*}

\elem

\bpf First recall that $|W_s^{*,\e_N}(\alpha)-V_s^{*,\e_N}({\bf {V}}_{s},  {\bf {W}}_{s},  \alpha)|^2
=|W_s^{*}(\alpha)-V_s^{*}({\bf {V}}_{s},  {\bf {W}}_{s},  \alpha)|^2$, see Notation \ref{nnn}.
It thus follows from \eqref{ee3} (with $v=W^1_s$, $v_*=W_s^{*,\e_N}(\alpha)$, $\tilde v=V^1_s$
and $\tilde v_*=V_s^{*,\e_N}({\bf {V}}_{s},  {\bf {W}}_{s},  \alpha)$) that
\begin{align*}
I_0^N(s) &\le C\int_0^1\Big(|W_s^1-V_s^{1}|^2+|W_s^{*}(\alpha)-V_s^{*}({\bf {V}}_{s},  {\bf {W}}_{s},  \alpha)|^2\Big)|W_s^1-W_s^{*,\e_N}(\alpha)|^\gamma d\alpha \\
&\hskip 5cm +C
K^{1-2/\nu}\int_0^1 |W_s^1-W_s^{*,\epsilon_N}(\alpha)|^{2+2\gamma/\nu} d\alpha.
\end{align*}

Next, we study $I_1 ^N(s)$. As seen in the proof of  Lemma \ref{estimategeneral},
$$
\int_0^\infty\!\!\int_0^{2\pi}\!\!c(v,v_*,z,\varphi) d\varphi dz=-(v-v_*)\Phi(|v-v_*|)
\;\hbox{ and }\;\int_0^\infty\!\!\int_0^{2\pi}\!\!c_K(v,v_*,z,\varphi) d\varphi dz=-(v-v_*)\Phi_K(|v-v_*|),
$$
where $\Phi(x)= \pi\int_0^{\infty}(1-\cos G(z/x^\gamma))dz$ and $\Phi_K(x)= \pi\int_0^{K}(1-\cos G(z/x^\gamma))dz$.
Then,
 \begin{align*}
I_1^N(s)
 &=2(W_s^1-V_s^{1}) \cdot \int_0^1\Big[-\big(W_{s}^{1}-W_s^*(\alpha)\big)\Phi\big(|W_{s}^{1}-W_s^*(\alpha)|\big)\\
 &\hskip3cm~ +\big(W_s^1-W_s^{*,\epsilon_N}(\alpha)\big)\Phi\big(|W_s^1-W_s^{*,\epsilon_N}(\alpha)|\big)\\
 &\hskip3cm~ -\big(V_s^{1}-V_s^{*,\epsilon_N}({\bf {V}}_{s},  {\bf {W}}_{s},  \alpha)\big)\Phi_{K}\big(|V_s^{1}-V_s^{*,\epsilon_N}({\bf {V}}_{s},  {\bf {W}}_{s},  \alpha)|\big)\\
 &\hskip3cm~ +\big(V_s^{1}-V_s^{*}({\bf {V}}_{s},  {\bf {W}}_{s},  \alpha)\big)\Phi_{K}\big(|V_s^{1}-V_s^{*}({\bf {V}}_{s},  {\bf {W}}_{s},  \alpha)|\big)\Big]d\alpha.
\end{align*}
But we have checked that
$\left|X\Phi_K(|X|)- Y\Phi_K(|Y|)\right| \leq C |X-Y||X|^\gamma$ for any $X,Y\in\bb{R}^3$
in the proof of  Lemma \ref{estimategeneral}, and it of course also holds true that
$\left|X\Phi(|X|)- Y\Phi(|Y|)\right| \leq C |X-Y||X|^\gamma$. Thus
\begin{align*}
I_1^N(s)\leq & C |W_s^1-V_s^{1}| \int_0^1
\Big[ |W_s^*(\alpha)-W_s^{*,\epsilon_N}(\alpha)||W_s^1-W_s^{*,\epsilon_N}(\alpha)|^\gamma\\
&\hskip3cm+|V_s^{*,\epsilon_N}({\bf {V}}_{s},  {\bf {W}}_{s},  \alpha)-V_s^{*}({\bf {V}}_{s},  {\bf {W}}_{s},  \alpha)||V_s^{1}-V_s^{*,\epsilon_N}({\bf {V}}_{s},  {\bf {W}}_{s},  \alpha)|^\gamma
\Big]d\alpha\\
=& C |W_s^1-V_s^{1}| \int_0^1 |\e_NY(\alpha)| \Big[|W_s^1-W_s^{*}(\alpha)-\e_NY(\alpha)|^\gamma
+|V_s^{1}-V_s^{*}({\bf {V}}_{s},  {\bf {W}}_{s},  \alpha)-\e_NY(\alpha)|^\gamma
\Big]d\alpha\\
\leq & C  |W_s^1-V_s^{1}|^2\\
&+C \e_N^2 \int_0^1 |Y(\alpha)|^2 \Big[|W_s^1-W_s^{*}(\alpha)-\e_NY(\alpha)|^{2\gamma}
+|V_s^{1}-V_s^{*}({\bf {V}}_{s},  {\bf {W}}_{s},  \alpha)-\e_NY(\alpha)|^{2\gamma}
\Big]d\alpha.
\end{align*}
But $Y$ is independent of $(W_s^{*},V_s^{*}({\bf {V}}_{s},  {\bf {W}}_{s},\cdot))$
and it holds that $\sup_{x \in \mathbb{R}^3} \int_0^1|x-\e_N Y(\alpha)|^{2\gamma}|Y(\alpha)|^2d\alpha\leq
\int_0^1|\e_N Y(\alpha)|^{2\gamma}|Y(\alpha)|^2d\alpha=C \e_N^{2\gamma}$
(recall that $\gamma \in (-1,0)$ and that $Y$ is uniformly distributed on $B(0,1)$), so that finally,
$$
I_1^N(s)\leq  C  |W_s^1-V_s^{1}|^2 + C \e_N^{2+2\gamma}.
$$

For $I_2^N (s)$, we first write $I_2^N (s)\leq A+B$, where
$$
A=2\int_0^1\int_0^\infty\int_0^{2\pi} |c_W(s)-c_W^{N}(s)|^2d\varphi dzd\alpha
\;\hbox{ and}\;B=2\int_0^1\int_0^\infty\int_0^{2\pi}|c_{K,V}^{N}(s)-c_{K,V}(s)|^2d\varphi dzd\alpha.
$$
We first apply \eqref{ee2} with with $v=W^1_s$, $v_*=W_s^{*,\e_N}(\alpha)$, $\tilde v=W^1_s$
and $\tilde v_*=W_s^{*}(\alpha)$:
$$
A\leq C \int_0^1 |W_s^*(\alpha)-W_s^{*,\epsilon_N}(\alpha)|^2 |W_s^{1}-W_s^{*,\epsilon_N}(\alpha)|^\gamma  d\alpha
=C \e_N^2  \int_0^1 |Y(\alpha)|^2|W_s^{1}-W_s^{*}(\alpha)-\e_NY(\alpha)|^\gamma  d\alpha.
$$
Using that
$\sup_{x \in \mathbb{R}^3} \int_0^1|x-\e_N Y(\alpha)|^{\gamma}|Y(\alpha)|^2d\alpha\leq
\int_0^1|\e_N Y(\alpha)|^{\gamma}|Y(\alpha)|^2d\alpha=C \e_N^{\gamma}$ and
arguing as in the study of $I^N_1(s)$, we conclude that
$A \leq C \e_N^{2+\gamma} \leq C \e_N^{2+2\gamma}$.
The other term $B$ is treated in the same way
(observe that \eqref{ee2} obviously also holds when replacing $c$ by $c_K=c{\bf 1}_{\{z\leq K\}}$).

\vskip1mm

We finally treat $I_3^N (s)$.  It is obvious that
\begin{align*}
I_3^N(s)&\le \int_0^1\int_0^\infty\int_0^{2\pi}  |c_W^N (s)-c_{K,V}^N (s)|^2 d\varphi dzd\alpha+I_2^N(s).
\end{align*}
But
\begin{align*}
\int_0^\infty\int_0^{2\pi}  |c_W^N (s)-c_{K,V}^N (s)|^2 d\varphi dz
=& \int_0^K\int_0^{2\pi}  |c_W^N(s)-c_{V}^N(s)|^2 d\varphi dz+ \int_K^\infty\int_0^{2\pi}  |c_W^N(s)|^2 d\varphi dz.
\end{align*}
Applying first \eqref{ee2} with  $v=W^1_s$, $v_*=W_s^{*,\e_N}(\alpha)$, $\tilde v=V^1_s$
and $\tilde v_*=V_s^{*,\e_N}({\bf {V}}_{s},  {\bf {W}}_{s},  \alpha)$, we find that
\begin{align*}
&\int_0^K\int_0^{2\pi}  |c_W^N(s)-c_{V}^N (s)|^2 d\varphi dz\\
\le& C\big(|W_s^1-V_s^{1}|^2+|W_s^{*,\epsilon_N}(\alpha)-V_s^{*,\epsilon_N}({\bf {V}}_{s},  {\bf {W}}_{s},  \alpha)|^2\big) |W_s^1-W_s^{*,\epsilon_N}(\alpha)|^\gamma\\
=& C\big(|W_s^1-V_s^{1}|^2+|W_s^{*}(\alpha)-V_s^{*}({\bf {V}}_{s},  {\bf {W}}_{s},  \alpha)|^2\big) |W_s^1-W_s^{*,\e_N}(\alpha)|^\gamma.
\end{align*}
Moreover, as seen in the proof of Lemma \ref{estimategeneral},
$\int_K^\infty\int_0^{2\pi}  |c_W^N (s)|^2 d\varphi dz=|W_s^1-W_s^{*,\epsilon_N}(\alpha)|^2
\Psi_K(|W_s^1-W_s^{*,\epsilon_N}(\alpha)|)$, where
$\Psi_K(x)=\Phi(x)-\Phi_K(x)\le C\int_K^\infty G^2(z/x^\gamma)dz\le Cx^{2\gamma/\nu}K^{1-2/\nu}$.  Hence,
$$
\int_K^\infty\int_0^{2\pi}  |c_W^N(s)|^2 d\varphi dz \leq C |W_s^1-W_s^{*,\epsilon_N}(\alpha)|^{2+2\gamma/\nu}K^{1-2/\nu}.
$$
All this shows that
\begin{align*}
I_3^N(s)\leq& I_2^N(s) + C\int_0^1\big(|W_s^1-V_s^{1}|^2+|W_s^{*}(\alpha)-V_s^{*}({\bf {V}}_{s},  {\bf {W}}_{s},  \alpha)|^2\big) |W_s^1-W_s^{*,\e_N}(\alpha)|^\gamma d\alpha\\
& + CK^{1-2/\nu}\int_0^1
|W_s^1-W_s^{*,\epsilon_N}(\alpha)|^{2+2\gamma/\nu} d\alpha
\end{align*}
and this ends the proof.
\epf

To prove our main result, we need the following estimate which can be found in \cite[Theorem 1]{MR3383341}.

\blem\label{lem-empirical}
Fix $A>0$ and $q>4$. There is a constant $C_{A,q}$ such that for all
$f\in\ca{P}_q(\bb{R}^3)$ verifying $\int_{\bb{R}^3}|v|^q f(dv)\le A$, all i.i.d. family
$(X_i)_{i=1,...,N}$ of $f$-distributed random variables,
\begin{align*}
\bb{E}\left[\ca{W}_2^2\Big(f, N^{-1}\sum_{i=1}^N\delta_{X_i}\Big)\right]\le  C_{A,q} N^{-1/2}.
\end{align*}
\elem

\bprop\label{pro-distance}
Fix $T>0$ and recall that  $h_N$ was defined in Lemma \emph{\ref{dicrete-norm}}. Consider the stopping time
\[\sigma_N=\inf\{t\ge0: \|\bar\mu_{{\bf {W}}_t}^N\|_{L^p}\ge 13500(1+h_N(t))\},\]
where $\bar\mu_{{\bf {W}}_t}^N=\mu_{{\bf {W}}_t}^N *\psi_{\e_N}$ with $\psi_{\e_N}(x)=(3/(4\pi\epsilon_N^3))\bbd{1}_{\{|x|\le \epsilon_N\}}$ and $\mu_{{\bf {W}}_t}^N=N^{-1}\sum_{1}^N \delta_{W^i_t}$. We have
for  all $T>0$,
\[\sup_{[0,T]}\bb{E}[|W_{t\wedge\sigma_N} ^1-V_{t\wedge\sigma_N} ^{1}|^2] \le C_{T}(\epsilon_N^{2+2\gamma}+K^{1-2/\nu}+N^{-1/2}).\]

\eprop

\bpf
We fix $T>0$ and set $u_t^{N}=\bb{E}[|W_{t\wedge\sigma_N} ^1-V_{t\wedge\sigma_N} ^{1}|^2]$ for all $t\in [0,T]$.
By the It\^{o} formula, we have
\begin{align*}
u_t^{N}
&=\bb{E}\Big[\int_0^{t\wedge\sigma_N} \int_0^1\int_0^\infty \int_0^{2\pi}\Big(|W_{s}^1-V_{s}^{1}+c_W(s)-c_{K,V}(s)|^2-|W_{s}^1-V_{s}^1|^2\Big)d\varphi dzd\alpha\Big]\\
&=\bb{E}\Big[\int_0^{t\wedge\sigma_N} \int_0^1\int_0^\infty \int_0^{2\pi}\Big(2(W_{s}^1-V_{s}^{1})\cdot(c_W(s)-c_{K,V}(s))+|c_W(s)-c_{K,V}(s)|^2\Big)d\varphi dzd\alpha\Big]\\
&=\bb{E}\Big[\int_0^{t\wedge\sigma_N} \big(I_0^{N}(s)+I_1^{N}(s)+I_2^{N}(s)+I_3^{N}(s)\big) ds\Big],
\end{align*}
where $I_i^{N}(s)$  was introduced in Lemma \ref{intermediate-lemma} for $i=0,1,2,3$.
We know from Lemma \ref{intermediate-lemma}  that
\begin{align*}
u_t ^{N} \le& Ct\epsilon_N^{2+2\gamma} + C\int_0^t u_s^N ds + C(J^N_1(t)+J^N_2(t)+J^N_3(t)),
\end{align*}
where
\begin{align*}
J^N_1(t)=&
\bb{E}\Big[\int_0^{t\wedge\sigma_N} \int_0^1 |W_s^1-V_s^{1}|^2|W_s^1-W_s^{*,\epsilon_N}(\alpha)|^\gamma d\alpha ds\Big],\\
J^N_2(t)=&\bb{E}\Big[\int_0^{t\wedge\sigma_N} \int_0^1|W_s^{*}(\alpha)-V_s^{*}({\bf {V}}_{s},  {\bf {W}}_{s},  \alpha)|^2|W_s^1-W_s^{*,\epsilon_N}(\alpha)|^\gamma d\alpha ds\Big],\\
J^N_3(t)=&K^{1-2/\nu}\bb{E}\Big[\int_0^{t\wedge\sigma_N} \int_0^1|W_s^1-W_s^{*,\epsilon_N}(\alpha)|^{2+2\gamma/\nu} d\alpha ds\Big].
\end{align*}

First, we have
$$
J^N_3(t)\leq C K^{1-2/\nu} t.
$$
Indeed, it suffices to use that $|W_s^1-W_s^{*,\epsilon_N}(\alpha)|^{2+2\gamma/\nu}
\leq C(1+|W^1_s|^2+|W_s^{*,\epsilon_N}(\alpha)|^2)$ (because $2+2\gamma/\nu \in (0,2)$),
that $|W_s^{*,\epsilon_N}(\alpha)|^2\leq 2+2|W_s^{*}(\alpha)|^2$ (because $\e_N\in (0,1)$
and $Y$ takes its values in $B(0,1)$) and finally that $\E[|W^1_s|^2]=\int_0^1 |W_s^{*}(\alpha)|^2d\alpha
=m_2(f_0)$.

\vskip1mm

Next, $\ca{L}_\alpha(W_s^{*,\e_N})=f_s*\psi_{\e_N}$, so that
$\int_0^1|W_s^1-W_s^{*,\epsilon_N}(\alpha)|^{\gamma} d\alpha \le 1+C_{\gamma,p}\|f_s*\psi_{\e_N}\|_{L^p}$
by  (\ref{norm-inequality}) (recall that $p>3/(3+\gamma)$ is fixed since the begining of the section).
Of course, $\|f_s*\psi_{\e_N}\|_{L^p}\le \|f_s\|_{L^p}$, and we conclude that
\begin{align*}
J^N_1(t)\leq
C_{\gamma,p} \int_0^{t} (1+\|f_s\|_{L^p}) \, u_s^{N}\  ds.
\end{align*}

On the other hand, using the exchangeability  and  that
$W_s^{*,\epsilon_N}(\alpha)=W_s^{*}(\alpha)+\epsilon_N Y(\alpha)$,   with $Y(\alpha)$
independent of $W_s^{*}(\alpha)$ and $V_s^{*}({\bf {V}}_{s},  {\bf {W}}_{s},  \alpha)$ introduced
in Notation \ref{nnn},
we have
\begin{align*}
J^N_2(t)
&=\bb{E}\Big[\int_0^{t\wedge\sigma_N} \int_0^1|W_s^{*}(\alpha)-V_s^{*}({\bf {V}}_{s},  {\bf {W}}_{s},  \alpha)|^2
N^{-1}\sum_{i=1}^N \Big|W_s^i-\epsilon_N Y(\alpha)-W_s^{*}(\alpha) \Big|^\gamma d\alpha ds\Big]\\
&=\bb{E}\Big[\int_0^{t\wedge\sigma_N} \int_0^1
|W_s^{*}(\alpha)-V_s^{*} ({\bf {V}}_{s},  {\bf {W}}_{s},  \alpha)|^2
\Big(\int_{\bb{R}^3}\int_{\bb{R}^3}|w-x-W_s^{*}(\alpha)|^\gamma \psi_{\epsilon_N}(x)\mu_{{\bf {W}}_s}^N (dw) dx \Big) d\alpha ds\Big]\\
&=\bb{E}\Big[\int_0^{t\wedge\sigma_N} \int_0^1
|W_s^{*}(\alpha)-V_s^{*}({\bf {V}}_{s},  {\bf {W}}_{s},  \alpha)|^2
\Big(\int_{\bb{R}^3}|w-W_s^{*}(\alpha)|^\gamma \bar\mu_{{\bf {W}}_s}^N(dw) \Big) d\alpha ds\Big].
\end{align*}
But $\int_{\bb{R}^3}|W_s^{*}(\alpha)-w|^\gamma \bar\mu_{{\bf {W}}_s}^N(dw) \le C_{\gamma,p} (1+\|\bar\mu_{{\bf {W}}_s}^N\|_{L^p})$
by \eqref{norm-inequality}, so that
\begin{align*}
J^N_2(t)
&\le C_{\gamma,p}\bb{E}\Big[\int_0^{t\wedge\sigma_N} \int_0^1(1+\|\bar\mu_{{\bf {W}}_s}^N\|_{L^p})|W_s^{*}(\alpha)-V_s^{*}({\bf {V}}_{s},  {\bf {W}}_{s},  \alpha)|^2 d\alpha ds \Big].
\end{align*}
We now deduce from Lemma \ref{coupling} that
\begin{align*}
\int_0^1 |W_s^{*}(\alpha)-V_s^{*}({\bf {V}}_{s},  {\bf {W}}_{s},  \alpha)|^2 d\alpha
 & \le 2
\int_0^1\Big(|W_s^{*}(\alpha)-Z_s^{*}({\bf {W}}_{s},  \alpha)|^2
+ |Z_s^{*}( {\bf {W}}_{s},  \alpha) - V_s^{*}({\bf {V}}_{s},  {\bf {W}}_{s},  \alpha)|^2\Big) d\alpha\\
& = 2 \ca{W}_2^2( f_s, \mu_{{\bf {W}}_s}^N ) + 2 \frac{1}{N}\sum_{i=1}^{N}|W_s^i - V_s^i|^2.
\end{align*}
Using the exchangeability and that $\|\bar\mu_{{\bf {W}}_s}^N\|_{L^p}\le 13500(1+h_N(s))$ for all $s\le \tau_N$,
it holds that
\begin{align*}
J^N_2(t)
&\le C \int_0^t(1+h_N(s))\bb{E}[\ca{W}_2^2( f_s, \mu_{{\bf {W}}_s}^N )]ds+ C \int_0^t(1+h_N(s))\, u_s^N ds.
\end{align*}

We thus have checked that
$$
u_t^{N}
 \le C t (\epsilon_N^{2+2\gamma} + K^{1-2/\nu} ) + C \int_0^t \big(1+h_N(s) \big)\bb{E}\big[\ca{W}_2^2( f_s, \mu_{{\bf {W}}_s}^N )\big]ds + C \int_0^t \big( 1+\|f_{s}\|_{L^p}+h_N(s) \big) u_s^{N}ds.
$$
But for each $t\geq 0$, the family
$(W^i_t)_{i=1,\dots,N}$ is i.i.d. and $f_t$-distributed. Furthermore,  $\sup_{[0,T]} \E[|W^1_t|^q]<\infty$
($q>6$) by \eqref{lem-momen}. Hence Lemma \ref{lem-empirical} tells  us that
\begin{equation}\label{ddd}
\sup_{[0,T]}\bb{E}\big[\ca{W}_2^2( f_s, \mu_{{\bf {W}}_s}^N )\big] \leq C_{T}N^{-1/2}.
\end{equation}
Using the Gr\"{o}nwall lemma, we deduce that
\begin{align*}
\sup_{[0,T]} u_t^{N}
& \le C_T\left(\epsilon_N^{2+2\gamma}+K^{1-2/\nu}+N^{-1/2} \int_0^T(1+h_N(s))ds \right)\exp{\Big(C\int_0^T (1+\|f_s\|_{L^p}+h_N(s))ds\Big)}.
\end{align*}
But  $\int_0^T h_N(s)ds \leq 2\int_0^T \|f_s\|_{L^p} ds$ by Lemma \ref{dicrete-norm}-(ii).
And we know that $f\in L^1_{loc}\big([0,\infty), L^{p}({\bb{R}^3})\big)$. We thus conclude that
\begin{align*}
\sup_{[0,T]} u_t^{N} & \le C_T\left(\epsilon_N^{2+2\gamma}+K^{1-2/\nu}+N^{-1/2}\right)
\end{align*}
as desired.
\epf

Now, we give the
\bpf[Proof of Theorem \ref{main-result}]
As explained at the beginning of the section, we only have to prove \eqref{ggooal}.
Recall that $\sigma_N=\inf\{t\ge0: \|\bar\mu_{{\bf {W}}_t}^N\|_{L^p}\ge 13500(1+h_N(t))\}$,
that $q>6$ and that $\delta=6/q$. It is clear that
$\bb{P}[\sigma_N\le T]\le C_{T,q,\delta}N^{1-q\delta/3}=C_{T,q}N^{-1}$ from Proposition \ref{norm-bound}.
Then for $t\in [0,T]$, we write
\begin{align*}
\sup_{[0, T]}\bb{E}[\ca{W}_2^2(\mu_{{\bf {V}}_t}^N, f_t)]
\le 2\sup_{[0, T]}\bb{E}[\ca{W}_2^2(\mu_{{\bf {V}}_t}^N, \mu_{{\bf {W}}_t}^N)+\ca{W}_2^2(\mu_{{\bf {W}}_t}^N, f_t)]
\leq 2\sup_{[0, T]}\bb{E}[\ca{W}_2^2(\mu_{{\bf {V}}_t}^N, \mu_{{\bf {W}}_t}^N)]+C_T N^{-1/2}
\end{align*}
by \eqref{ddd}.
But, by exchangeability, we have
$$
\bb{E}[\ca{W}_2^2(\mu_{{\bf {V}}_t}^N, \mu_{{\bf {W}}_t}^N)]
\le \bb{E}\Big[N^{-1}\sum_{i=1}^N|W_{t} ^i-V_{t} ^{i}|^2\Big]=\bb{E}[|W_{t} ^1-V_{t} ^{1}|^2].
$$
Moreover,
\begin{align*}
\bb{E}[|W_{t} ^1-V_{t} ^{1}|^2]
&\le \bb{E}[|W_{t\wedge\sigma_N} ^1-V_{t\wedge\sigma_N} ^{1}|^2] +\bb{E}[|W_{t} ^1-V_{t} ^{1}|^2\bbd{1}_{\{\sigma_N \le T\}}]\\
&\le C_T(\epsilon_N^{2+2\gamma}+K^{1-2/\nu}+N^{-1/2})+ C\bb{E}[|W_{t} ^1|^4+|V_{t} ^{1}|^4]^{1/2}(\bb{P}(\sigma_N \le T))^{1/2},
\end{align*}
by Proposition \ref{pro-distance} , and the Cauchy-Schwarz inequality.
Noting that  $\bb{E}[|W_{t} ^1|^4]\le C_T$ by \eqref{lem-momen}, and that $\bb{E}[|V_{t} ^{1}|^4]\le C_T \bb{E}[|V_0^1|^4]$ by Lemma \ref{aaaa}, we deduce that
$$
\bb{E}[|W_{t} ^1-V_{t} ^{1}|^2]
\le C_{T,q}(\epsilon_N^{2+2\gamma}+K^{1-2/\nu}+N^{-1/2}).
$$
All in all, we have proved that
\begin{align*}
\sup_{[0, T]}\bb{E}[\ca{W}_2^2(\mu_{{\bf {V}}_t}^N, f_t)]\leq C_{T,q}(\epsilon_N^{2+2\gamma}+K^{1-2/\nu}+N^{-1/2}).
\end{align*}
This is precisely \eqref{ggooal}, since $\epsilon_N^{2+2\gamma}=N^{-(1-6/q)(2+2\gamma)/3}$, with
$\e_N=N^{-(1-\delta)/3}$ and $\delta=6/q$.
\epf

\section*{Acknowledgements}
\quad~~
I would like to thank greatly Nicolas Fournier for continuous and generous supports in this research and especially Maxime Hauray for inspiring discussion about the proof of Proposition 3.1-Step 1.

\bibliographystyle{abbrv} %plain / abbrv /amsalpha/siam
\bibliography{psb}
\vskip 3mm
\textsc{Liping XU\\
Laboratoire de Probabilit\'es et Mod\`eles Al\'eatoires, UMR 7599, Universit\'e Pierre-et-Marie Curie,
Case 188, 4 place Jussieu, F-75252 Paris Cedex 5, France.}\\
\emph{E-mail}: \texttt{liping.xu@upmc.fr}.

\end{document}